\documentclass[a4paper,11pt]{article}                                                   

\usepackage[english]{babel}
\usepackage{amsfonts,amsmath,amsthm,graphicx,a4wide}   

\usepackage[colorlinks=false]{hyperref}
\usepackage{hyperref}
\usepackage{graphicx} 
\usepackage{amsmath} 
\usepackage{amssymb}  
\usepackage{mathrsfs}
\usepackage{bbm}
\usepackage{xcolor}
\usepackage{mathrsfs}
\usepackage{dsfont} 
\usepackage{authblk}

\newtheorem{theorem}{Theorem}[section]
\newtheorem{lemma}[theorem]{Lemma}
\newtheorem{proposition}[theorem]{Proposition} 

\newtheorem{remark}[theorem]{Remark}
\newtheorem{definition}[theorem]{Definition}

\title{\LARGE \bf
Robust feedback stabilization of $N$-level quantum spin systems
}
\date{\vspace{-5ex}}

\author[1]{Weichao Liang}
\author[2]{Nina H. Amini}
\author[2]{Paolo Mason}

\affil[1]{Laboratoire Analyse G\'eom\'etrie Mod\'elisation, CY Cergy Paris Universit\'e, 2, av. Adolphe Chauvin, 95302 Cergy-Pontoise, cedex, France.
  {\tt\small weichao.liang@u-cergy.fr}.}
\affil[2]{Laboratoire des signaux et syst\`{e}mes (L2S), CNRS-CentraleSup\'{e}lec-Universit\'{e} Paris-Sud, Universit\'{e} Paris-Saclay, 3, rue Joliot Curie, 91190 Gif-sur-Yvette, France
  {\tt\small [first name].[family name]@l2s.centralesupelec.fr}.\footnote{\bf Funding: \rm This work is supported by Agence Nationale de la Recherche projects Q-COAST
ANR-19-CE48-0003 and QUACO ANR-17-CE40-0007.}}

\begin{document}

\maketitle
\begin{abstract}
In this paper, we consider $N$-level quantum angular momentum systems interacting with electromagnetic fields undergoing continuous-time measurements. We suppose unawareness of the initial state and physical parameters, entailing the introduction of an additional state representing the estimated quantum state. The evolution of the quantum state and its estimation  is described by a coupled stochastic master equation. Here, we study the asymptotic behavior of such a system in presence of a feedback controller. We provide sufficient conditions on the feedback controller and on the estimated parameters that guarantee exponential stabilization of the coupled stochastic system towards an eigenstate of the measurement operator. Furthermore, we estimate  the corresponding rate of convergence.  We also provide parametrized feedback laws satisfying such conditions. Our results show the robustness of the feedback stabilization strategy considered in~\cite{liang2019exponential} in case of imprecise initialization of the estimated state and with respect to the unknown physical parameters. 

\end{abstract}
\section{Introduction}
The evolution and feedback control of an open quantum system which undergoes continuous-time measurements can be studied in the framework of quantum stochastic calculus and quantum probability theory. These essential mathematical tools have been introduced by Hudson and Parthasarathy \cite{hudson1984quantum} in the 1980s. 
A primitive theory of quantum filtering theory was developed by Davies in the 1960s~\cite{davies1969quantum,davies1976quantum}. In the 1980s, Belavkin established the quantum filtering theory~\cite{belavkin1983theory,belavkin1989nondemolition,belavkin1992quantum,belavkin1995quantum} which is a more natural extension of the classical theory (see, e.g.,~\cite{kallianpur2013stochastic}). 
In this context, the quantum state, called quantum filter, represents the stochastic evolution of the conditional density operator when the
system interacts with the field. The theory developed by Belavkin provides an essential foundation for statistical inference in, for instance, quantum optical systems. 
In the physics community, a heuristic approach to quantum filtering, called quantum trajectory theory, has been developed by Carmichael in the early 1990s~\cite{carmichael1993open}.
For a modern review of quantum filtering based on the  Hudson-Parthasarathy quantum stochastic calculus, see \cite{bouten2007introduction}.

The measurement-based feedback as a branch of stochastic control was developed by Belavkin in \cite{belavkin1983theory}. Later, Bouten and van Handel established a separation principle \cite{bouten2008separation}, showing that, in order to design a state-based feedback, the quantum filtering problem and the control problem can be studied separately. This is an important contribution which provides the analogue of the important separation principle from standard optimal stochastic control. 

The evolution of an open quantum system undergoing indirect continuous-time measurements is described by the so-called quantum stochastic master equation. The deterministic part of this equation 
is given by the well known Lindblad operator. Its stochastic part represents the back-action effect of continuous-time measurements. 
For controlled open quantum systems,  the control inputs usually appear in the Lindblad operator (through the system Hamiltonian).
Feedback control of open quantum systems, in order to prepare pure states, has been the subject of many papers, e.g., \cite{van2005feedback, mirrahimi2007stabilizing, tsumura2008global, liang2019exponential}. 
The preparation of pure states is investigated as an essential step towards quantum technologies~\cite{sayrin2011real,hacker2019deterministic}.
The problem of finding a quantum feedback controller that globally stabilizes a quantum spin-$\frac 12$ system towards an eigenstate of the measurement operator $\sigma_z$ in the
presence of imperfect measurements was first tackled in~\cite{van2005feedback}. This feedback controller was designed by looking numerically for an appropriate global Lyapunov function. Later, in~\cite{mirrahimi2007stabilizing}, by analyzing the stochastic flow and using stochastic Lyapunov techniques, the authors constructed a switching control law to stabilize globally $N$-dimensional quantum angular momentum systems around any predetermined eigenstate of the measurement operator. Recently, in~\cite{liang2018exponential,liang2019exponential} we established exponential stabilization results for spin-$\frac 12$ and spin-$J$ systems  towards any stationary state of the open-loop dynamics by means of a continuous feedback. Our approach combined local stochastic stability analysis with the use of the support theorem. Unlike earlier works, based on the LaSalle approach, our techniques allowed us to estimate the rate of convergence to the target state~\cite{liang2018exponential,liang2019exponential}. This is important in view of practical implementation in quantum information processing. In~\cite{cardona2018exponential,cardona2020exponential}, the authors provided exponential stabilization results via a different approach.

 In real experiments, different types of imperfections, such as detection inefficiencies and unawareness of initial states, may be present (see e.g.,~\cite{sayrin2011real}). The design of stabilizing feedback controllers robust to such experimental imperfections is a crucial step towards engineering of quantum devices. 
In the case of unawareness of initial states, 
one considers an additional state representing an estimation of the actual quantum state. When a feedback is applied, this feedback depends on the estimated state and the observation process depends on the actual quantum state.
 This leads to a coupled stochastic master equation where, for the dynamics of the estimated state, one chooses an arbitrary initial state and fixes some estimated physical parameters. 
An important  question is whether the estimated process approaches the true quantum trajectory when time tends to infinity; such a property is often referred to as 
stability of quantum filter~\cite{van2009stability,amini2014stability,benoist2014large}. 

In~\cite{liang2019exponential_TwoUnknown}, we showed the convergence of the estimated state towards the true state for spin-$\frac 12$ systems under appropriate assumptions on the feedback and we conjectured the possibility  of stabilizing exponentially $N$-level quantum angular momentum systems by means of a candidate feedback in presence of unawareness of initial states and detection inefficiencies. 
 Recently, in~\cite{liang2020robustness_two}, still for spin-$\frac 12$ systems, we showed that the feedback strategy proposed in~\cite{liang2018exponential} robustly stabilizes the system in case of unawareness of the initial state and physical parameters.

In this paper, we study the feedback exponential stabilization problem for $N$-level quantum angular momentum systems  in the case of unknown initial states and imprecise knowledge of the physical parameters (the detection efficiency, the free Hamiltonian, and the strength of the interaction between the system and the probe). We provide sufficient conditions on the feedback controller and the estimated parameters ensuring exponential stabilization of the quantum state and its estimation towards a predetermined eigenstate of the measurement operator $J_z.$ Compared to previous literature and in particular~\cite{liang2019exponential}, where the initial states (and physical parameters) are known, the main difficulty here is the presence of additional equilibria for the coupled system which cannot be eliminated. Hence in order to prove the stabilization towards the target equilibrium, one needs to prove the instability of such additional equilibria and develop a reachability analysis for the coupled system starting outside these equilibria. To our knowledge, this study provides the first results on asymptotic and exponential stabilization of $N$-level quantum angular momentum systems  in the case of unknown initial states and imprecise knowledge of the physical parameters.

The paper is structured as follows. We first provide some preliminary results 
concerning invariance properties for the coupled system
(Section~\ref{sec:pre}). Secondly, we present general Lyapunov-type conditions ensuring almost sure exponential stabilization and asymptotic stabilization of the coupled system (Section~\ref{sec:general}). Based on the general results, we obtain explicit conditions on the estimated parameters and the feedback controller, guaranteeing almost sure exponential convergence and further providing an estimate for the corresponding convergence rate (Section~\ref{sec:explicit}). Finally, we design a parametrized family of feedback controllers satisfying such conditions (Section~\ref{sec:par}). This proves in particular~\cite[Conjecture~4.4]{liang2019exponential_TwoUnknown} in the case in which the target state corresponds to the first or the last eigenstate of the measurement operator $J_z,$ even in the case of unawareness of the physical parameters. Numerical simulations are provided in order to illustrate our results and to support the efficiency of the proposed candidate feedback (Section~\ref{sec:sim}).
\paragraph{Notations}
The imaginary unit is denoted by $i$. We take $\mathds{1}$ as the indicator function. Given a complex number $z$, we indicate as $\mathbf{Re}(z)$ its real part. We denote the conjugate transpose of a matrix $A$ by $A^*.$ The function $\mathrm{Tr}(A)$ corresponds to the trace of a square matrix $A.$ The commutator of two square matrices $A$ and $B$ is denoted by $[A,B]:=AB-BA.$ We denote by $\mathrm{int}(\mathcal{S})$ the interior of a subset of a topological space and by $\partial \mathcal{S}$ its boundary.
\section{System description and problem setting}
Here, we consider a $N$-level quantum spin system under continuous-time homodyne measurements with $N>1$. The stochastic master equations describing the evolution of the  
quantum state and the corresponding estimation
are given as follows,
\begin{align}
d\rho_t&=L^u_{\omega,M}(\rho_t)dt+G_{\eta,M}(\rho_t)\big(dY_t-2\sqrt{\eta M}\mathrm{Tr}(J_z\rho_t)dt\big), \label{Eq:ND SME}\\
d\hat{\rho}_t&=L^u_{\hat{\omega},\hat{M}}(\hat{\rho}_t)dt+G_{\hat{\eta},\hat{M}}(\hat{\rho}_t)\big(dY_t-2\sqrt{\hat{\eta} \hat{M}}\mathrm{Tr}(J_z\hat{\rho}_t)dt\big),\label{Eq:ND SME filter}
\end{align}
where
\begin{itemize}
\item the actual state of the $N$-level quantum spin system, denoted by $\rho$, belongs to the compact space
$
\mathcal{S}_N:=\{\rho\in\mathbb{C}^{N\times N}|\,\rho=\rho^*,\mathrm{Tr}(\rho)=1,\rho \geq 0\}.
$
The associated estimated state is denoted by $\hat{\rho}\in\mathcal{S}_N$,
\item the drift terms are given by $L^u_{\omega,M}(\rho):=-i[\omega J_z+uJ_y,\rho]+\frac{M}{2}(2J_z\rho J_z-J_z^2\rho-\rho J_z^2)$ and the diffusion terms are given by $G_{\eta,M}(\rho):=\sqrt{\eta M}\big(J_z\rho+\rho J_z-2\mathrm{Tr}(J_z\rho)\rho\big)$,
\item $Y_t$ denotes the observation process of the actual quantum system, which is a continuous semi-martingale with $\langle Y,Y\rangle_t=t$. Its dynamics satisfy $dY_t=dW_t+2\sqrt{\eta M}\mathrm{Tr}(J_z\rho_t)dt$, where $W_t$ is a one-dimensional standard Wiener process, 
\item $u:=u(\hat{\rho})\in\mathbb{R}$ denotes the continuous feedback controller, as a function of the estimated state $\hat{\rho}$, adapted to $\mathcal{F}^y_t:=\sigma(Y_s,0\leq s\leq t)$, which is the $\sigma$-field generated by the observation process up to time $t$,
\item 
$J_z$ is the (self-adjoint) angular momentum along the axis $z$, and it is defined by 
\begin{equation*}
J_z e_n=(J-n)e_n,\quad n\in\{0,\dots,2J\}, 
\end{equation*}
where $J:=\frac{N-1}{2}$ represents the fixed angular momentum and $\{e_0,\dots,e_{2J}\}$ corresponds to an orthonormal basis of $\mathbb C^N.$ With respect to this basis, the matrix form of $J_z$ is given by
\begin{equation*}
J_z=
\begin{bmatrix}
J &&&&  \\
& J-1&&&\\
&&\ddots&&\\
&&&-J+1&\\
&&&&-J
\end{bmatrix},
\end{equation*}
\item 
$J_y$ is the (self-adjoint) angular momentum along the axis $y$, and it is defined by 
\begin{equation*}
\begin{split}
J_ye_n=-ic_{n}e_{n-1}+ic_{n+1}e_{n+1},\quad n\in\{0,\dots,2J\},
\end{split}
\end{equation*}
where $c_m=\frac12\sqrt{(2J+1-m)m}$. The matrix form of $J_y$ is given by
\begin{equation*}
J_y=
\begin{bmatrix}
0&-ic_1 &&&\\
ic_1&0&-ic_2&&\\
&\ddots&\ddots&\ddots&\\
&&ic_{2J-1}&0&-ic_{2J}\\
&&&ic_{2J}&0
\end{bmatrix},
\end{equation*}  
\item $\eta\in(0,1]$ describes the efficiency of the detectors,  $M>0$ is the strength of the interaction between the system and the probe, and $\omega \geq 0$ is a parameter characterizing the free Hamiltonian.  We assume that these parameters are not precisely known in practice and that the estimated parameters are given by $\hat{\eta}\in(0,1]$, $\hat{M}>0$ and $\hat{\omega}\geq 0$. 
\end{itemize}
By replacing $dY_t=dW_t+2\sqrt{\eta M}\mathrm{Tr}(J_z\rho_t)dt$ in equations~\eqref{Eq:ND SME}--\eqref{Eq:ND SME filter}, we obtain the following matrix-valued stochastic differential equation in It\^o form, which describes the time evolution of the pair $(\rho_t,\hat{\rho}_t)\in\mathcal{S}_N\times\mathcal{S}_N$,
\begin{align}
d\rho_t=&L^u_{\omega,M}(\rho_t)dt+G_{\eta,M}(\rho_t)dW_t,\label{Eq:ND SME W}\\
d\hat{\rho}_t=&L^u_{\hat{\omega},\hat{M}}(\hat{\rho}_t)dt+G_{\hat{\eta},\hat{M}}(\hat{\rho}_t)dW_t\label{Eq:ND SME filter W}\\
&+2\sqrt{\hat{\eta} \hat{M}}G_{\hat{\eta},\hat{M}}(\hat{\rho}_t)\big(\sqrt{\eta M}\mathrm{Tr}(J_z\rho_t)-\sqrt{\hat{\eta} \hat{M}}\mathrm{Tr}(J_z\hat{\rho}_t)\big)dt.\nonumber
\end{align}
If $u\in\mathcal{C}^1(\mathcal{S}_N,\mathbb{R})$,  the existence and uniqueness of the solution of the coupled system~\eqref{Eq:ND SME W}--\eqref{Eq:ND SME filter W} can be shown by similar arguments as in~\cite[Proposition 3.5]{mirrahimi2007stabilizing}. Moreover, it can be shown as in~\cite[Proposition 3.7]{mirrahimi2007stabilizing} that $(\rho_t,\hat{\rho}_t)$ is a strong Markov process in $\mathcal{S}_N\times\mathcal{S}_N$.

Obviously, if $\rho_0=\hat{\rho}_0$, $\omega=\hat{\omega}$, $\eta=\hat{\eta}$ and $M=\hat{M}$, then $\rho_t=\hat{\rho}_t$ for all $t\geq 0$ almost surely. In this case, the state feedback stabilization of the system~\eqref{Eq:ND SME W} towards the target state $\boldsymbol \rho_{\bar{n}}:=e_{\bar{n}}e^*_{\bar{n}}$ with $\bar{n}\in\{0,\dots,2J\}$ has been studied in several papers (see e.g.,~\cite{mirrahimi2007stabilizing,liang2019exponential}). In particular, sufficient conditions on the state feedback controller guaranteeing the almost sure exponential stabilization has been provided in~\cite{liang2019exponential}. 
In this paper, we will study the following more general problem.

\medskip

\paragraph{Problem} Assume that we do not have access to the initial state $\rho_0$ and that the actual physical parameters $\omega,\eta,M$ are not precisely known. Find conditions on the estimated parameters $\hat{\omega},\hat{\eta},\hat{M}$ and the feedback controller $u(\hat{\rho}),$ which ensure the exponential convergence of the solutions of~\eqref{Eq:ND SME W}--\eqref{Eq:ND SME filter W} towards the target state $(\boldsymbol \rho_{\bar{n}},\boldsymbol \rho_{\bar{n}})$, with $\bar{n}\in\{0,\dots,2J\}$, independently of the choice of $\hat{\rho}_0$.
\section{Some basic tools for stochastic processes and stability}
In this section, we first recall some elementary tools for the study of stochastic processes, we then present some notions of stochastic stability and we finally state the support theorem~\cite{stroock1972support} which will be used throughout the paper.
\paragraph{Infinitesimal generator and It\^o formula}
Given a stochastic differential equation $dq_t=f(q_t)dt+g(q_t)dW_t$, where $q_t$ takes values in $Q\subset \mathbb{R}^p,$ the infinitesimal generator  is the operator $\mathscr{L}$ acting on twice continuously differentiable functions $V: Q \times \mathbb{R}_+ \rightarrow \mathbb{R}$ in the following way
\begin{equation*}
\mathscr{L}V(q,t):=\frac{\partial V(q,t)}{\partial t}+\sum_{i=1}^p\frac{\partial V(q,t)}{\partial q_i}f_i(q)+\frac12 \sum_{i,j=1}^p\frac{\partial^2 V(q,t)}{\partial q_i\partial q_j}g_i(q)g_j(q).
\end{equation*}
It\^o formula describes the variation of the function $V$ along solutions of the stochastic differential equation and is given as follows
\begin{equation*}
dV(q,t) = \mathscr{L}V(q,t)dt+\sum_{i=1}^p\frac{\partial V(q,t)}{\partial q_i}g_i(q)dW_t.
\end{equation*}
From now on, the operator $\mathscr{L}$ is associated with Equations~\eqref{Eq:ND SME W}--\eqref{Eq:ND SME filter W}.
 
\paragraph{Stochastic stability} We recall that the Bures distance~\cite{bengtsson2017geometry} between two density matrices $\rho^{(1)}$ and $\rho^{(2)}$ is given by
\begin{equation*}
d_B(\rho^{(1)},\rho^{(2)}):=\sqrt{2-2\sqrt{\mathcal{F}(\rho^{(1)},\rho^{(2)})}},
\end{equation*}
where $\mathcal{F}(\hat{\rho},\rho):=\mathrm{Tr}\big(\sqrt{\sqrt{\rho^{(1)}}\rho^{(2)}\sqrt{\rho^{(1)}}}\big)\in[0,1]$ is so-called the fidelity~\cite{nielsen2002quantum}.  In particular, the Bures distance between $\rho\in\mathcal{S}_N$ and a pure state $\boldsymbol\rho_{n}:=e_ne^*_n$ with $n\in\{0,\dots,2J\}$ is given by 
$
d_B(\rho,\boldsymbol \rho_n)=\sqrt{2-2\sqrt{\rho_{n,n}}},
$
where $\rho_{n,n}:=\mathrm{Tr}(\rho\boldsymbol\rho_{n})$ denotes the projection of $\rho$ on $\boldsymbol\rho_{n}$. We then define the neighbourhood $B_r(\rho)$ of $\rho\in\mathcal S_N$ as
$
B_r(\rho) := \{\sigma \in \mathcal{S}_N|\, d_B(\rho,\sigma) < r\}.
$

The distance between two elements in $\mathcal{S}_N\times\mathcal{S}_N$, which will be needed to adapt classical notions of stochastic stability  (see e.g.,~\cite{mao2007stochastic, khasminskii2011stochastic}) to our setting, is then defined as
\begin{equation*}
\mathbf{d}_B\big((\rho^{(1)},\hat{\rho}^{(1)}),(\rho^{(2)},\hat{\rho}^{(2)})\big):=d_B(\rho^{(1)},\rho^{(2)})+d_B(\hat{\rho}^{(1)},\hat{\rho}^{(2)}).
\end{equation*}
We denote the ball of radius $r$ around $(\rho,\hat{\rho})$ as
\begin{equation*}
\mathbf{B}_r(\rho,\hat{\rho}) := \{(\sigma,\hat{\sigma}) \in \mathcal{S}_N\times\mathcal{S}_N|\, \mathbf{d}_B\big((\rho,\hat{\rho}),(\sigma,\hat{\sigma})\big) < r\}.
\end{equation*}

\begin{definition}[see e.g., \cite{khasminskii2011stochastic,mao2007stochastic}]
Let $(\boldsymbol \rho,\hat{\boldsymbol \rho})$ be an equilibrium of the coupled system~\eqref{Eq:ND SME W}--\eqref{Eq:ND SME filter W}, then $(\boldsymbol \rho,\hat{\boldsymbol \rho})$ is said to be
\begin{enumerate}
\item[1.] 
\emph{locally stable in probability}, if for every $\varepsilon \in (0,1)$ and for every $r >0$, there exists  $\delta = \delta(\varepsilon,r)$ such that,
\begin{equation*}
\mathbb{P} \big( (\rho_t,\hat{\rho}_t) \in \mathbf{B}_r(\boldsymbol \rho,\hat{\boldsymbol \rho}) \text{ for } t \geq 0 \big) \geq 1-\varepsilon,
\end{equation*}
whenever $(\rho_0,\hat{\rho}_0) \in \mathbf{B}_{\delta}(\boldsymbol \rho,\hat{\boldsymbol \rho})$.

\item[2.]
\emph{almost surely asymptotically stable} in $\Gamma$, where $\Gamma\subset\mathcal{S}_N\times\mathcal{S}_N$ is a.s. invariant, if it is locally stable in probability and,
\begin{equation*}
\mathbb{P} \left( \lim_{t\rightarrow\infty}\mathbf{d}_B\big((\rho,\hat{\rho}),(\boldsymbol \rho,\hat{\boldsymbol \rho})\big)=0 \right) = 1,
\end{equation*}
whenever $(\rho_0,\hat\rho_0) \in \Gamma.$


\item[3.]
\emph{almost surely exponentially stable} in $\Gamma$, where $\Gamma\subset\mathcal{S}_N\times\mathcal{S}_N$ is a.s. invariant, if
\begin{equation*}
\limsup_{t \rightarrow \infty} \frac{1}{t} \log \mathbf{d}_B\big((\rho,\hat{\rho}),(\boldsymbol \rho,\hat{\boldsymbol \rho})\big) < 0, \quad a.s.
\end{equation*}
whenever $(\rho_0,\hat{\rho}_0) \in \Gamma.$ The left-hand side of the above inequality is called the \emph{sample Lyapunov exponent} of the solution.
\end{enumerate}
\end{definition}

\paragraph{Stratonovich equation and Support theorem}
Any stochastic differential equation in It\^o form in $\mathbb R^K$
\begin{equation*}
dx_t=\widehat X_0(x_t)dt+\sum^n_{k=1}\widehat X_k(x_t)dW^k_t, \quad x_0 = x,
\end{equation*}
can be written in the following Stratonovich form~\cite{rogers2000diffusions2}
\begin{equation*}
dx_t = X_0(x_t)dt+\sum^n_{k=1}X_k(x_t) \circ dW^k_t, \quad x_0 = x,
\end{equation*}
where 
$X_0(x)=\widehat X_0(x)-\frac{1}{2}\sum^K_{l=1}\sum^n_{k=1}\frac{\partial \widehat X_k}{\partial x_l}(x)(\widehat X_k)_l(x)$, $(\widehat X_k)_l$ denoting the component $l$ of the vector $\widehat X_k,$ and $X_k(x)=\widehat X_k(x)$ for $k\neq 0$.

\medskip
The following classical theorem relates the solutions of a stochastic differential equation with those of an associated deterministic one.
\begin{theorem}[Support theorem~\cite{stroock1972support}]
Let $X_0(t,x)$ be a bounded measurable function, uniformly Lipschitz continuous in $x$ and $X_k(t,x)$  be continuously differentiable in $t$ and twice continuously differentiable in $x$, with bounded derivatives, for $k\neq 0.$ Consider the Stratonovich equation
\begin{equation*}
dx_t = X_0(t,x_t)dt+\sum^n_{k=1}X_k(t,x_t) \circ dW^k_t, \quad x_0 = x,
\end{equation*}
and denote by $\mathbb{P}_x$ the probability law of the solution $x_t$ starting at $x$. Consider in addition the associated deterministic control system
\begin{equation*}
\frac{d}{dt}x_{v}(t) = X_0(t,x_{v}(t))+\sum^n_{k=1}X_k(t,x_{v}(t))v^k(t), \quad x_v(0) = x,
\end{equation*}
with $v^k \in \mathcal{V}$, where $\mathcal{V}$ is the set of all locally bounded measurable functions from $\mathbb{R}_+$ to $\mathbb{R}$. Define $\mathcal{W}_x$ as the set of all continuous paths from $\mathbb{R}_+$ to $\mathbb R^K$ starting at $x$, equipped with the topology of uniform convergence on compact sets, and $\mathcal{I}_x$ as the smallest closed subset of $\mathcal{W}_x$ such that $\mathbb{P}_x(x_{\cdot} \in \mathcal{I}_x)=1$. Then, $\mathcal{I}_x = \overline{ \{ x_{v}(\cdot)\in\mathcal{W}_x|\, v \in \mathcal{V}^n\} } \subset \mathcal{W}_x.$
\label{Thm:Support}
\end{theorem}
\section{Feedback stabilization of the coupled system}
In this section, we provide conditions on the feedback controller $u(\hat{\rho})$ and a suitable domain of the estimated parameters $\hat{\omega}$, $\hat{M}$ and $\hat{\eta}$, which ensure the exponential stabilization of $(\rho_t,\hat{\rho}_t)$ towards a target state $(\boldsymbol \rho_{\bar{n}},\boldsymbol \rho_{\bar{n}})$ with $\bar{n}\in\{0,\dots,2J\}$.


We impose the following hypothesis, which implies that the coupled system~\eqref{Eq:ND SME W}--\eqref{Eq:ND SME filter W} contains exactly the $N$ equilibria $(\boldsymbol\rho_n,\boldsymbol\rho_{\bar{n}})$ with $n\in\{0,\dots,2J\}$. 
\medskip

\begin{itemize}
\item[\textbf{H0}:] 
$u\in\mathcal{C}^1(\mathcal{S}_N,\mathbb{R})$, $u(\boldsymbol\rho_{\bar{n}})=0$ and $u(\hat{\rho})\neq0$ for all $\hat{\rho}\in\{\boldsymbol\rho_0,\dots,\boldsymbol\rho_{2J}\}\setminus \boldsymbol\rho_{\bar{n}}$.
\end{itemize}

 \medskip

\begin{remark}
It can be easily verified that, if we turn off the feedback controller, i.e., $u(\hat{\rho})\equiv0$, there are $N^2$ equilibria $(\boldsymbol\rho_n,\boldsymbol\rho_m)$ with $n,m\in\{0,\dots,2J\}$ for the coupled system~\eqref{Eq:ND SME W}--\eqref{Eq:ND SME filter W}. In~\cite[Theorem 7]{benoist2014large}, the authors show that, for the case $u\equiv0$, $\eta=\hat{\eta}=1$, $\omega=\hat{\omega}$ and $M=\hat{M}$, the trajectories of the coupled system~\eqref{Eq:ND SME W}--\eqref{Eq:ND SME filter W} converge exponentially almost surely towards the subset  $\{(\boldsymbol\rho_0,\boldsymbol\rho_0),\dots,(\boldsymbol\rho_{2J},\boldsymbol\rho_{2J})\}$ of the set of the equilibria.
\end{remark}

In the following, we first provide some preliminary technical results. Secondly, we present general Lyapunov-type conditions ensuring exponential stabilization of the coupled system towards $(\boldsymbol\rho_{\bar{n}},\boldsymbol\rho_{\bar{n}}).$ We then apply this result to obtain an explicit stabilization result and an estimated rate of exponential convergence, under suitable assumptions on the estimated parameters $\hat\eta,\hat M$ and on the feedback controller $u(\hat{\rho})$. Finally, we design a parametrized family of feedback controllers satisfying such conditions.
\subsection{Preliminary results}~\label{sec:pre}
Before starting our analysis, we state some fundamental results that will be needed later. These results are analogous to the results in~\cite[Section~4]{liang2019exponential} and they concern invariance properties  for the coupled system~\eqref{Eq:ND SME W}--\eqref{Eq:ND SME filter W}  involving the boundary 
$\partial\mathcal{S}_N=\{\rho\in\mathcal{S}_N|\,\det({\rho})=0\}$ and the interior  $\mathrm{int}(\mathcal{S}_N)=\{\rho\in\mathcal{S}_N|\,{\rho}>0\}.$ Since their proofs are based on the same arguments, we omit them. 
\begin{lemma}
Assume that $u\in\mathcal{C}^1(\mathcal{S}_N,\mathbb{R})$. Let $(\rho_t,\hat{\rho}_t)$ be a solution of the coupled system~\eqref{Eq:ND SME W}--\eqref{Eq:ND SME filter W}. If $\rho_0\in\mathrm{int}(\mathcal{S}_N)$, then 
$\mathbb{P}(\rho_t\in\mathrm{int}(\mathcal{S}_N),\,\forall t\geq0)=1.$ 
Similarly, if $\hat{\rho}_0\in\mathrm{int}(\mathcal{S}_N)$, then $\mathbb{P}(\hat{\rho}_t\in\mathrm{int}(\mathcal{S}_N),\,\forall t\geq0)=1.$ 
More in general, the ranks of ${\rho}_t$ and $\hat{\rho}_t$ are a.s. non-decreasing.
\label{Lemma:PosDef invariant}
\end{lemma}
\begin{lemma}
Assume that $u\in\mathcal{C}^1(\mathcal{S}_N,\mathbb{R})$. 
If $\eta=1$, then $\partial \mathcal{S}_N\times \mathcal{S}_N$ is a.s. invariant for the coupled system~\eqref{Eq:ND SME W}--\eqref{Eq:ND SME filter W}. If $\hat{\eta}=1$, then $\mathcal{S}_N\times\partial\mathcal{S}_N$ is a.s. invariant for the coupled system~\eqref{Eq:ND SME W}--\eqref{Eq:ND SME filter W}.
\label{Lemma:Boundary invariant}
\end{lemma}
\begin{lemma}
Assume that \emph{\textbf{H0}} is satisfied. Let $(\rho_t,\hat{\rho}_t)$ be a solution of the coupled system~\eqref{Eq:ND SME W}--\eqref{Eq:ND SME filter W}. If 
$\hat{\rho}_0 \neq \boldsymbol \rho_{\bar{n}}$, then $\mathbb{P}( \hat{\rho}_t \neq \boldsymbol \rho_{\bar{n}}, \forall\, t\geq 0 )=1.$
Moreover, if $\hat{\rho}_0=\boldsymbol\rho_{\bar{n}}$ and 
$\rho_0\neq\boldsymbol \rho_n$, then $\mathbb{P}( \rho_t \neq \boldsymbol \rho_n, \forall\, t\geq 0 )=1,$ for any $n\in\{0,\dots,2J\}$. 
\label{Lemma:Never Reach Lemma}
\end{lemma}
\subsection{General results on exponential stabilization and asymptotic stabilization}\label{sec:general}
In order to obtain our general results, we first provide sufficient conditions on the feedback controller guaranteeing the exponential instability of the equilibria $(\boldsymbol\rho_n,\boldsymbol\rho_{\bar{n}})$ with $n\neq\bar{n}$. Secondly, we show that, under suitable conditions on the feedback controller, for all initial states, except $(\boldsymbol\rho_n,\boldsymbol\rho_{\bar{n}})$ for $n\neq\bar{n}$, the trajectories $(\rho_t,\hat{\rho}_t)$ can enter in an arbitrary neighbourhood of the target state $(\boldsymbol\rho_{\bar{n}},\boldsymbol\rho_{\bar{n}})$ in finite time almost surely (reachability property). Thirdly, we establish a general result ensuring almost sure exponential convergence under some assumptions on the feedback controller and additional local Lyapunov-type conditions. An analogous result concerning  almost sure asymptotic stabilization is also provided under weaker 
Lyapunov-type conditions.
\subsubsection{Exponential instability of the equilibria $(\boldsymbol\rho_n,\boldsymbol\rho_{\bar{n}})$ with $\bf n\neq\bar{n}$} Analogous to the instability results for classical stochastic system~\cite[Theorem 4.3.5]{mao2007stochastic}, in the following lemmas, we show the exponential instability of the equilibria $(\boldsymbol\rho_n,\boldsymbol\rho_{\bar{n}})$, with $n\neq\bar{n}$, separately for the cases $\bar{n}\in\{0,2J\}$ and $\bar{n}\in\{1,\dots,2J-1\}$. In order to state such results, we need to introduce two assumptions.

\medskip
\begin{itemize}
\item[\textbf{H1}:] 
$|u(\hat{\rho})|\leq c(1-\hat{\rho}_{\bar{n},\bar{n}})^m$ for $\hat\rho\in\mathcal S_N$ with $m>1/2$ and for some constant $c>0$.
\end{itemize}
\medskip

The above assumption is required to show the instability of $(\boldsymbol\rho_n,\boldsymbol\rho_{\bar{n}})$ with $n\neq\bar{n}$ and $\bar{n}\in\{0,2J\}.$ We remark that the above hypothesis is not a consequence of the continuous differentiability of $u.$
\medskip

\begin{itemize}
\item[\textbf{H2}:] 
$u(\hat{\rho})=0$ for all $\hat{\rho}\in B_{\xi}(\boldsymbol\rho_{\bar{n}})$, for a sufficiently small $\xi>0$.
\end{itemize}
\medskip

The above assumption is needed to show the instability of $(\boldsymbol\rho_n,\boldsymbol\rho_{\bar{n}})$ with $n\neq\bar{n}$ and the reachability of  an arbitrary neighbourhood of the target state $(\boldsymbol\rho_{\bar{n}},\boldsymbol\rho_{\bar{n}})$ for  $\bar{n}\notin\{0,2J\}.$ 
\medskip

Furthermore, for $n\in\{0,\dots,2J\},$ we set $\Theta_{{n}}(\rho):=\mathrm{Tr}(i[J_y,\rho]\boldsymbol\rho_{{n}})$, $P_{{n}}(\rho):=J-{n}-\mathrm{Tr}(J_z\rho)$, and 
$\mathcal{T}(\rho,\hat{\rho}):=\sqrt{\eta M}\mathrm{Tr}(J_z\rho)-\sqrt{\hat{\eta}\hat{M}}\mathrm{Tr}(J_z\hat{\rho}).$
\begin{lemma}
Let $\bar{n}\in\{0,2J\}$. Assume that $\hat{\rho}_0\in\mathcal{S}_N\setminus\{\boldsymbol\rho_{\bar{n}}\}$,  \emph{\textbf{H0}} and  \emph{\textbf{H1}} are satisfied  and 
\begin{equation}
(N-2)\sqrt{\hat{\eta}\hat{M}}>(N-3)\sqrt{\eta M}.
\label{eq:robustness-0,2J}
\end{equation}
Then, there exists $r>0$ such that, for all $(\rho_0,\hat{\rho}_0)\in\mathbf{B}_r(\boldsymbol\rho_n,\boldsymbol\rho_{\bar{n}})$ with $n\neq \bar{n}$, the trajectories of the coupled system~\eqref{Eq:ND SME W}--\eqref{Eq:ND SME filter W} exit $\mathbf{B}_r(\boldsymbol\rho_n,\boldsymbol\rho_{\bar{n}})$ in finite time almost surely.
\label{Lemma:Instability 0 and 2J}
\end{lemma}
\proof
Consider the function $V_{\bar{n}}(\hat{\rho})=1-\hat{\rho}_{\bar{n},\bar{n}}\in[0,1]$.
By \textbf{H1}, we have that $|u|\leq c V_{\bar{n}}(\hat{\rho})^m$ with $m>1/2$ and for some constant $c>0$. Moreover 
\begin{align*}
|\Theta_{\bar{n}}(\hat{\rho})| & = |2c_{\bar n+1} \mathbf{Re}(\hat{\rho}_{\bar n,\bar n+1})-2c_{\bar n} \mathbf{Re}(\hat{\rho}_{\bar n,\bar n-1})|\\
&\leq 2(c_{\bar n}+c_{\bar n+1}) V_{\bar{n}}(\hat{\rho})^{1/2},
\end{align*}
with the convention that $c_0=c_{N}=0$.
It is also easy to see that $P_0(\hat{\rho})\geq V_0(\hat{\rho})$ and $P_{2J}(\hat{\rho})\leq -V_{2J}(\hat{\rho})$.
Under the assumption 
\begin{equation}
J\sqrt{\hat{\eta}\hat{M}}>(J-1)\sqrt{\eta M}
\label{eq:asp}
\end{equation}
 we have that $\mathcal{T}(\rho,\hat{\rho})<0$ for $\bar{n}=0$, and $\mathcal{T}(\rho,\hat{\rho})>0$ for $\bar{n}=2J$, in some small enough neighbourhood of $(\boldsymbol\rho_n,\boldsymbol\rho_{\bar{n}})$ with $n\neq \bar{n}$.
Hence, in such a neighbourhood,
\begin{equation*}
\begin{split}
\mathscr{L}V_{\bar{n}}(\hat{\rho})&=u\Theta_{\bar{n}}(\hat{\rho})-4\sqrt{\hat{\eta}\hat{M}}\mathcal{T}(\rho,\hat{\rho})P_{\bar{n}}(\hat{\rho})\hat{\rho}_{\bar{n},\bar{n}}\\
&\geq \big(-\alpha V_{\bar{n}}(\hat{\rho})^{m-\frac12}+4\sqrt{\hat{\eta}\hat{M}}\left|\mathcal{T}(\rho,\hat{\rho})\right|\hat{\rho}_{\bar{n},\bar{n}}\big)V_{\bar{n}}(\hat{\rho}),
\end{split}
\end{equation*}
with some $\alpha>0$. This implies $\liminf_{(\rho,\hat{\rho})\rightarrow(\boldsymbol\rho_n,\boldsymbol\rho_{\bar{n}})}\frac{\mathscr{L}V_{\bar{n}}(\hat{\rho})}{V_{\bar{n}}(\hat{\rho})}\geq C_1$, where 

\noindent $C_1:=4\sqrt{\hat{\eta}\hat{M}}\left( J\sqrt{\hat{\eta}\hat{M}}-(J-1)\sqrt{\eta M}\right).$ 
We have
\begin{equation*}
\mathscr{L}\log V_{\bar{n}}(\hat{\rho})=\frac{\mathscr{L}V_{\bar{n}}(\hat{\rho})}{V_{\bar{n}}(\hat{\rho})}-\frac{1}{2}\left(\frac{\partial V_{\bar{n}}(\hat{\rho})}{\partial \hat{\rho}}\frac{G_{\hat{\eta},\hat{M}}(\hat{\rho})}{V_{\bar{n}}(\hat{\rho})}\right)^2.
\end{equation*}
We want to guarantee that the right-hand side of the previous expression is larger than a positive constant around the equilibrium $(\boldsymbol\rho_n,\boldsymbol\rho_{\bar{n}})$. For this purpose we notice that 
$$\limsup_{(\rho,\hat{\rho})\rightarrow(\boldsymbol\rho_n,\boldsymbol\rho_{\bar{n}})}\left(\frac{\partial V_{\bar{n}}(\hat{\rho})}{\partial \hat{\rho}}\frac{G_{\hat{\eta},\hat{M}}(\hat{\rho})}{V_{\bar{n}}(\hat{\rho})}\right)^2\leq C_2,$$ where $C_2:=4\hat{\eta}\hat{M}.$ The condition $(2J-1)\sqrt{\hat{\eta}\hat{M}}>(2J-2)\sqrt{\eta M}$ guarantees that $C_1>C_2/2$. This condition implies~\eqref{eq:asp} and, by replacing $J=\frac{N-1}{2}$, it reduces to $(N-2)\sqrt{\hat{\eta}\hat{M}}>(N-3)\sqrt{\eta M}$. Note that, for $N\in\{2,3\}$, this condition is satisfied automatically. The relation $C_1-C_2/2>0$ implies that there exist $r>0$ and $C_3>0$ such that $\mathscr{L}\log V_{\bar{n}}\geq C_3$ for all $(\rho,\hat{\rho})\in \mathbf{B}_r(\boldsymbol\rho_n,\boldsymbol\rho_{\bar{n}})$. 

Let $\tau:=\inf\{t\geq0|\,(\rho_t,\hat{\rho}_t)\notin \mathbf{B}_r(\boldsymbol\rho_n,\boldsymbol\rho_{\bar{n}})\}$. Due to \hyperref[Lemma:Never Reach Lemma]{Lemma~\ref*{Lemma:Never Reach Lemma}}, we can apply It\^o's formula on $\log V_{\bar{n}}(\hat{\rho}_t)$. By taking the expectation\footnote{
Here and in the following, $\mathbb{P}_0$ corresponds to the probability law of $(\rho_t,\hat{\rho}_t)$ starting at $(\rho_0,\hat{\rho}_0)$; the associated expectation is denoted by $\mathbb{E}_0$}, we obtain the following  
\begin{equation*}
0\geq \mathbb{E}_0(\log V_{\bar{n}}(\hat{\rho}_{\tau}))\geq \log V_{\bar{n}}(\hat{\rho}_{0})+C_3\mathbb{E}_0(\tau).
\end{equation*}
Since $\hat{\rho}_0\neq\boldsymbol{\rho}_{\bar{n}}$, $\log V_{\bar{n}}(\hat{\rho}_{0})\in(-\infty,0]$. Thus, 
$
\mathbb{E}_0(\tau)\leq -\log V_{\bar{n}}(\hat{\rho}_{0})/C_3<\infty.
$
Then by Markov inequality, for all $(\rho_0,\hat{\rho}_0)\in\mathbf{B}_r(\boldsymbol\rho_n,\boldsymbol\rho_{\bar{n}})$ with $\hat{\rho}_0\neq\boldsymbol\rho_{\bar{n}}$, we have 
\begin{equation*}
\mathbb{P}_0(\tau=\infty) = \lim_{l\rightarrow \infty} \mathbb{P}_{0}(\tau \geq l) \leq \lim_{l\rightarrow \infty} \mathbb{E}_{0}(\tau)/l=0,
\end{equation*}
which implies
$
\mathbb{P}_{0}( \tau_{r}<\infty )=1.
$
The proof is complete.\hfill$\square$

Now, we show the 
instability of the equilibria $(\boldsymbol\rho_n,\boldsymbol\rho_{\bar{n}})$ for $n\neq \bar{n}$ and $\bar n$ not necessarily equal to $\{0,2J\}$.
\begin{lemma}
Let $\bar{n}\in\{0,\dots,2J\}$. Assume that $\hat{\rho}_0>0,$ \emph{\textbf{H0}} and \emph{\textbf{H2}} are satisfied  and 
\begin{equation}
(N-1)\sqrt{\eta M}>(N-2)\sqrt{\hat{\eta}\hat{M}}>(N-3)\sqrt{\eta M}.
\label{eq:parametersr}
\end{equation}
Then, there exists $r>0$ such that, for all $(\rho_0,\hat{\rho}_0)\in\mathbf{B}_r(\boldsymbol\rho_n,\boldsymbol\rho_{\bar{n}})$ with $n\neq \bar{n}$, the trajectories of the coupled system~\eqref{Eq:ND SME W}--\eqref{Eq:ND SME filter W} exit $\mathbf{B}_r(\boldsymbol\rho_n,\boldsymbol\rho_{\bar{n}})$ in finite time almost surely.
\label{Lemma:Instability General}
\end{lemma}
\proof
By \hyperref[Lemma:PosDef invariant]{Lemma~\ref*{Lemma:PosDef invariant}}, $\hat{\rho}_0>0$ implies $\hat{\rho}_t>0$ for all $t\geq0$, and therefore $(\hat{\rho}_t)_{k,k}>0$ for all $k\in\{0,\dots,2J\}$.
Consider the function $V_{\bar{n}}(\hat{\rho})=\hat{\rho}_{n,n}\in(0,1]$ with $n\neq \bar{n}$, whose infinitesimal generator is given by $\mathscr{L}V_{\bar{n}}(\hat{\rho})=-u\Theta_{n}(\hat{\rho})+4\sqrt{\hat{\eta}\hat{M}}P_{n}(\hat{\rho})\mathcal{T}(\rho,\hat{\rho})V_{\bar{n}}(\hat{\rho}).$
Note that, when $(\rho,\hat{\rho})$ converges to $(\boldsymbol\rho_n,\boldsymbol\rho_{\bar{n}})$, we have that $P_{n}(\hat{\rho}),\mathrm{Tr}(J_z\rho),\mathrm{Tr}(J_z\hat{\rho})$ converge to $\bar{n}-n,J-n,J-\bar n$, respectively.
By \textbf{H2} we thus get the following estimate
\begin{equation*}
\liminf_{(\rho,\hat{\rho})\rightarrow(\boldsymbol\rho_n,\boldsymbol\rho_{\bar{n}})}\frac{\mathscr{L}V_{\bar{n}}(\hat{\rho})}{V_{\bar{n}}(\hat{\rho})}\geq C_1,\quad C_1:= 4\sqrt{\hat{\eta}\hat{M}}(\bar{n}-n)\left((J-n)\sqrt{\eta M} - (J-\bar{n})\sqrt{\hat{\eta}\hat{M}}\right),
\end{equation*}
on a small enough neighbourhood of $(\boldsymbol\rho_n,\boldsymbol\rho_{\bar{n}})$.
Moreover, we have 
\[\limsup_{(\rho,\hat{\rho})\rightarrow(\boldsymbol\rho_n,\boldsymbol\rho_{\bar{n}})}\left(\frac{\partial V_{\bar{n}}(\hat{\rho})}{\partial \hat{\rho}}\frac{G_{\hat{\eta},\hat{M}}(\hat{\rho})}{V_{\bar{n}}(\hat{\rho})}\right)^2\leq C_2,\qquad C_2:=4\hat{\eta}\hat{M}(\bar{n}-n)^2.\]
A sufficient condition on $\hat{\eta}\hat{M}$ and $\eta M$ to guarantee that $C_1-C_2/2>0$ is $(2J-2)\sqrt{\eta M}<(2J-1)\sqrt{\hat{\eta}\hat{M}}<2J\sqrt{\eta M}$. By replacing $J=\frac{N-1}{2}$, we obtain $(N-3)\sqrt{\eta M}<(N-2)\sqrt{\hat{\eta}\hat{M}}<(N-1)\sqrt{\eta M}$. The relation $C_1-C_2/2>0$ implies that there exist $r>0$ small enough and $C_3>0$ such that, for all $(\rho,\hat{\rho})\in \mathbf{B}_r(\boldsymbol\rho_n,\boldsymbol\rho_{\bar{n}})$,
\begin{equation*}
\mathscr{L}\log V_{\bar{n}}(\hat{\rho})=\frac{\mathscr{L}V_{\bar{n}}(\hat{\rho})}{V_{\bar{n}}(\hat{\rho})}-\frac{1}{2}\left(\frac{\partial V_{\bar{n}}(\hat{\rho})}{\partial \hat{\rho}}\frac{G_{\hat{\eta},\hat{M}}(\hat{\rho})}{V_{\bar{n}}(\hat{\rho})}\right)^2>C_3>0.
\end{equation*}
The rest of the proof follows the same arguments as in Lemma~\ref{Lemma:Instability 0 and 2J}.\hfill$\square$
\begin{remark}
In the proof of~\hyperref[Lemma:Instability 0 and 2J]{Lemma~\ref*{Lemma:Instability 0 and 2J}}  it is shown that $\hat{\rho}_t$ almost surely exits a neighbourhood of $\boldsymbol\rho_{\bar{n}}$ by directly inspecting the variation of $(\hat{\rho}_t)_{\bar{n},\bar{n}}$. This method, however, can be applied only in the case $\bar{n}\in\{0,2J\}$ since otherwise the sign of $P_{\bar{n}}(\hat{\rho})$ (and therefore the sign of the infinitesimal generator of $\hat{\rho}_{\bar{n},\bar{n}}$) is not constant on a neighbourhood of the equilibrium.
In contrast, in~\hyperref[Lemma:Instability General]{Lemma~\ref*{Lemma:Instability General}}, dealing with the general case $\bar{n}\in\{0,\dots,2J\}$, we prove the instability of $(\boldsymbol\rho_n,\boldsymbol\rho_{\bar{n}})$ by showing that $(\hat{\rho}_t)_{n,n}$  moves away from zero, which indirectly implies that $(\hat{\rho}_t)_{\bar n,\bar n}$ moves away from one, when $(\rho_0,\hat{\rho}_0)$ approaches $(\boldsymbol\rho_n,\boldsymbol\rho_{\bar{n}})$. 
Not surprisingly, compared to~\hyperref[Lemma:Instability 0 and 2J]{Lemma~\ref*{Lemma:Instability 0 and 2J}},  \hyperref[Lemma:Instability General]{Lemma~\ref*{Lemma:Instability General}} needs stronger assumptions.
\end{remark}
\subsubsection{Reachability for the deterministic coupled system}
Define the following coupled deterministic system corresponding to the Stratonovich form of the coupled stochastic system~\eqref{Eq:ND SME W}--\eqref{Eq:ND SME filter W},
\begin{align}
\dot{\rho}_v(t)&=\tilde{L}^u_{\omega,\eta,M}\big(\rho_v(t)\big)+G_{\eta,M}\big(\rho_v(t)\big)V(t),\label{Eq:ND ODE}\\
\dot{\hat{\rho}}_v(t)&=\tilde{L}^u_{\hat{\omega},\hat{\eta},\hat{M}}\big(\hat{\rho}_v(t)\big)+G_{\hat{\eta},\hat{M}}\big(\hat{\rho}_v(t)\big)V(t),\label{Eq:ND ODE filter}
\end{align}
with $\rho_v(0)=\rho_0$, $\hat{\rho}_v(0)=\hat{\rho}_0$,
$\tilde{L}^u_{\omega,\eta,M}(\rho):=-i[\omega J_z+uJ_y,\rho]+M\big((1-\eta)J_z\rho J_z-\frac{1+\eta}{2}(J_z^2\rho+\rho J_z^2)+2\eta\mathrm{Tr}(J_z^2\rho)\rho\big)$, $V(t):=v(t)+2\sqrt{\eta M}\mathrm{Tr}\big(J_z\rho_v(t)\big)$ and $v(t)\in\mathcal{V}$ is the bounded control input.

Inspired by~\cite[Lemma 6.1]{liang2019exponential} and~\cite[Theorem 4.7]{baxendale1991invariant}, in the following lemmas, we analyze the possibility of constructing trajectories of~\eqref{Eq:ND ODE}--\eqref{Eq:ND ODE filter}  which enter 
an arbitrarily small neighbourhood of the target state.

Before stating the results, we define $\widehat{\mathbf{P}}_{{n}}:=\{\hat{\rho}\in\mathcal{S}_N|\,P_{{n}}(\hat{\rho})=0\}$ for $n\in\{0,\dots,2J\}$ and the ``variance function'' $\mathscr{V}_z(\hat{\rho}):=\mathrm{Tr}(J^2_z\hat{\rho})-\mathrm{Tr}(J_z\hat{\rho})^2$ of $J_z$ for the estimated state.
\begin{lemma}
Assume that $\bar{n}\in\{0,2J\}$ and  \emph{\textbf{H0}} holds true.   In addition,  assume that for any $\hat{\rho}_0\in\{\hat\rho\in\mathcal S_N| \, \hat{\rho}_{\bar n,\bar n}=0\},$  there exists a control $v(t)\in\mathcal{V}$ such that for all $t\in(0,\delta),$ with $\delta>0$ sufficiently small, $u\big(\hat{\rho}_{v}(t)\big)\neq 0$, for some solution $\hat{\rho}_{v}(t)$ of Equation~\eqref{Eq:ND ODE filter}.
Then, for all $\varepsilon>0$ and any given initial state $(\rho_0,\hat{\rho}_0) \in\left(\mathcal{S}_{N}\times\mathcal{S}_{N}\right)\setminus\left\{\mathbf{B}_{\varepsilon}(\boldsymbol\rho_{\bar{n}},\boldsymbol\rho_{\bar{n}})\cup\bigcup_{n\neq \bar{n}}(\boldsymbol\rho_{n},\boldsymbol\rho_{\bar{n}})\right\},$ there exist $T\in(0,\infty)$ and $v(t)\in\mathcal{V}$ such that the trajectory $(\rho_v(t),\hat{\rho}_{v}(t))$ of the coupled deterministic system~\eqref{Eq:ND ODE}--\eqref{Eq:ND ODE filter} enters $\mathbf{B}_{\varepsilon}(\boldsymbol\rho_{\bar n},\boldsymbol\rho_{\bar n})$ for $t<T$.
\label{Lemma:Reachability ODE - 1}
\end{lemma}
\proof
From~\eqref{Eq:ND ODE}--\eqref{Eq:ND ODE filter} we have
\begin{small}
\begin{align}
\big(\dot{\rho}_v(t)\big)_{\bar{n},\bar{n}}=&-u\Theta_{\bar{n}}\big(\rho_v(t)\big)+2\sqrt{\eta M}\Lambda_{\bar{n}}\big(\rho_v(t)\big)\big(\rho_v(t)\big)_{\bar{n},\bar{n}}\label{Eq:ODE target}\\
&+4\sqrt{\eta M}P_{\bar{n}}\big(\rho_v(t)\big)\big(\rho_v(t)\big)_{\bar{n},\bar{n}}V(t),\nonumber\\
\big(\dot{\hat{\rho}}_v(t)\big)_{\bar{n},\bar{n}}=&-u\Theta_{\bar{n}}\big(\hat{\rho}_v(t)\big)+2\sqrt{\hat{\eta} \hat{M}}\Lambda_{\bar{n}}\big(\hat{\rho}_v(t)\big)\big(\hat{\rho}_v(t)\big)_{\bar{n},\bar{n}}\label{Eq:ODE  target filter}\\
&+4\sqrt{\hat{\eta} \hat{M}}P_{\bar{n}}\big(\hat{\rho}_v(t)\big)\big(\hat{\rho}_v(t)\big)_{\bar{n},\bar{n}}V(t),\nonumber
\end{align}
\end{small}
where $\rho_v(0)=\rho_0$, $\hat{\rho}_v(0)=\hat{\rho}_0$ and $\Lambda_{\bar{n}}(\rho):=\mathrm{Tr}(J^2_z\rho)-(J-\bar{n})^2$.
If $(\rho_0)_{\bar n,\bar n}=0$, by following the arguments of~\cite[Proposition 4.5]{liang2019exponential}, one can show the existence of a control input $v\in\mathcal{V}$ such that $(\hat{\rho}_{v}(t))_{\bar{n},\bar{n}}>0$ and $(\rho_{v}(t))_{\bar{n},\bar{n}}>0$ for $t\in(0,\delta)$. Thus, without loss of generality, we suppose $(\hat{\rho}_0)_{\bar{n},\bar{n}}>0$ and $(\rho_0)_{\bar{n},\bar{n}}>0.$ Moreover, we have $\widehat{\mathbf{P}}_{\bar{n}}= \boldsymbol\rho_{\bar n}$. Since $\mathcal{S}_{N}\times\mathcal{S}_{N}$ is compact, the first two terms of the right-hand side of Equations~\eqref{Eq:ODE target}--\eqref{Eq:ODE  target filter} are bounded from above in this domain and, as $\bar n\in\{0,2J\}$, one has $|P_{\bar n}(\rho)|\geq |1-\rho_{\bar n,\bar n}|>0$. Then by choosing $V=K/\min{(P_{\bar{n}}({\rho}),P_{\bar{n}}(\hat{\rho}))}$, with $K>0$ sufficiently large, we can guarantee that $(\rho_{v}(t),\hat{\rho}_{v}(t))\in \mathbf{B}_{\varepsilon}(\boldsymbol\rho_{\bar n},\boldsymbol\rho_{\bar n})$ for $t< T$ with $T<\infty$ if $(\hat{\rho}_0)_{\bar n,\bar n}>0$ and $(\rho_0)_{\bar n,\bar n}>0$. 
\hfill$\square$
\begin{lemma}
Assume that $\bar{n}\in\{1,\dots,2J-1\}$, \emph{\textbf{H0}} and \emph{\textbf{H2}} are satisfied. Let the parameters $\eta, M, \hat\eta, \hat M$ satisfy   
\begin{equation}
(N-1)\sqrt{\eta M}>(N-3)\sqrt{\hat{\eta}\hat{M}},\quad (N-1)\sqrt{\hat{\eta}\hat{M}}>(N-3)\sqrt{\eta M}.
\label{eq:parametersrf}
\end{equation}
Suppose, moreover, that
\begin{equation}
\forall \hat{\rho}\in\widehat{\mathbf{P}}_{\bar n}\setminus \boldsymbol\rho_{\bar n},\quad 2\hat{\eta} \hat{M}\mathscr{V}_{z}(\hat{\rho})\hat{\rho}_{\bar n,\bar n}>u(\hat{\rho})\Theta_{\bar{n}}(\hat{\rho}).
\label{Condition u}
\end{equation}
Then, for all $\varepsilon>0$ and any given initial state $(\rho_0,\hat{\rho}_0) \in (\mathcal{S}_N\times\mathrm{int}(\mathcal{S}_N))\setminus\mathbf{B}_{\varepsilon}(\boldsymbol\rho_{\bar{n}},\boldsymbol\rho_{\bar{n}})$, there exist $T\in(0,\infty)$ and $v(t)\in\mathcal{V}$ such that the trajectory $(\rho_v(t),\hat{\rho}_{v}(t))$ of the coupled deterministic system~\eqref{Eq:ND ODE}--\eqref{Eq:ND ODE filter} enters $\mathbf{B}_{\varepsilon}(\boldsymbol\rho_{\bar n},\boldsymbol\rho_{\bar n})$ with $t<T$.
\label{Lemma:Reachability ODE - 2}
\end{lemma}
\proof
Without loss of generality we assume in the following that $\varepsilon<\xi,$ where $\xi$ is as in $\textbf{H2}$. The proof proceeds in three steps:
\begin{enumerate}
\item First, we show that, for all $\hat{\rho}_0>0$,  there exists  $v\in\mathcal{V}$ such that $u(\hat{\rho}_v(t))\neq 0$ for some $t>0$.
\item Next, we show that, there exists $T_1\in(0,\infty)$ and $v\in\mathcal{V}$ such that $\hat{\rho}_v(T_1)\in B_{\varepsilon}(\boldsymbol\rho_{\bar{n}})$.
\item Finally, we show that, for any $\varepsilon\in(0,\xi)$, there exist $T_2\in(0,\infty)$ and $v(t)\in\mathcal{V}$ such that, $(\rho_v(t),\hat{\rho}_v(t))$ enters $\mathbf{B}_{\varepsilon}(\boldsymbol\rho_{\bar{n}},\boldsymbol\rho_{\bar{n}})$ with $t<T_2$.
\end{enumerate}

\emph{Step 1}: Suppose that $u(\hat{\rho}_v(t))= 0$ for every $v\in\mathcal{V}$ and $t\geq 0$.
Then, for any $n\neq\bar{n}$, we have  
\begin{equation*}
\big(\dot{\hat{\rho}}_v(t)\big)_{n,n}=2\sqrt{\hat{\eta} \hat{M}}\Lambda_{n}\big(\hat{\rho}_v(t)\big)\big(\hat{\rho}_v(t)\big)_{n,n}+4\sqrt{\hat{\eta} \hat{M}}P_{n}\big(\hat{\rho}_v(t)\big)\big(\hat{\rho}_v(t)\big)_{n,n}V(t).
\end{equation*}
By considering the two cases $\hat{\rho}\in\widehat{\mathbf{P}}_n\setminus B_{\varepsilon}(\boldsymbol\rho_n)$ with $\varepsilon>0$ sufficiently small and $\hat{\rho}\in\mathcal{S}_N\setminus\widehat{\mathbf{P}}_n$, and by employing the same argument as in the proof of~\cite[Lemma 6.1]{liang2019exponential}, we can show that there exists $v(t)\in\mathcal{V}$ such that $\hat{\rho}(t)$ enters $B_{\varepsilon}(\boldsymbol\rho_n)$. Moreover, $u(\hat{\rho})\neq0$ for all $\hat{\rho}\in B_{\varepsilon}(\boldsymbol\rho_n)$ due to \textbf{H0}, which leads to a contradiction. At once $u(\hat{\rho}_{t_0})\neq 0$ and, 
by the same arguments as in the proof of~\cite[Proposition 4.5]{liang2019exponential},
there exists a control input $v\in\mathcal{V}$ such that $\rho_{v}(t)>0$ for all $t>t_0$.

\emph{Step 2}: Due to \hyperref[Lemma:PosDef invariant]{Lemma~\ref*{Lemma:PosDef invariant}}, $\hat{\rho}_0>0$ implies $\hat{\rho}_t>0$ for all $t\geq 0$, thus $(\hat{\rho}_{t})_{\bar{n},\bar{n}}>0$ for all $t\geq0$. By the \hyperref[Thm:Support]{support theorem~\ref*{Thm:Support}}, we have $(\hat{\rho}_{v}(t))_{\bar{n},\bar{n}}>0$ for all $t\geq0$ and $v\in\mathcal{V}$. Then, by employing the same argument as in the proof of~\cite[Lemma 6.1]{liang2019exponential}, we can show that, for any $\varepsilon\in(0,\xi)$, there exists $T_1\in(0,\infty)$ and $v\in\mathcal{V}$ such that $\hat{\rho}_v(T_1)\in B_{\varepsilon}(\boldsymbol\rho_{\bar{n}})$. 

\emph{Step 3}: In this step, we show that for $\varepsilon\in(0,\xi)$ and $T_1\in(0,\infty)$ fixed in \emph{Step~2}, there exists $v\in\mathcal{V}$ such that $\hat{\rho}_v(t)$ remains in $B_{\varepsilon}(\boldsymbol\rho_{\bar{n}})$ for all $t\geq T_1$, while $\rho_v(t)$ enters $B_{\varepsilon}(\boldsymbol\rho_{\bar{n}})$. Due to the results established in \emph{Step~1} and \emph{Step~2}, we can assume $(\rho_v(t))_{\bar{n},\bar{n}}>0$ and $(\hat{\rho}_v(t))_{\bar{n},\bar{n}}>0$ for all $t\geq T_1$. Then, we define $f_n(t):=\frac{(\rho_v(t))_{n,n}}{(\rho_v(t))_{\bar{n},\bar{n}}}$ and $\hat{f}_n(t):=\frac{(\hat{\rho}_v(t))_{n,n}}{(\hat{\rho}_v(t))_{\bar{n},\bar{n}}}$. For all $n\neq \bar{n}$, we have
\begin{align}
\dot{\hat{f}}_n(t)=\frac{u(\hat{\rho}_v(t))}{(\hat{\rho}_v(t))_{\bar{n},\bar{n}}^2}&\big(\Theta_{\bar{n}}(\hat{\rho}_v(t))(\hat{\rho}_v(t))_{n,n}-\Theta_{n}(\hat{\rho}_v(t))(\hat{\rho}_v(t))_{\bar{n},\bar{n}}\big)\label{Eq:ODE f_hat}\\
&+2\sqrt{\hat{\eta} \hat{M}}(\bar{n}-n)\big(\sqrt{\hat{\eta} \hat{M}}(\bar{n}+n-2J)+2V(t)\big)\hat{f}_n(t).\nonumber
\end{align}
If $V(t)$ satisfies the following inequality
\begin{equation}
\sqrt{\hat{\eta} \hat{M}}(J-\bar{n}-\kappa)<V(t)<\sqrt{\hat{\eta} \hat{M}}(J-\bar{n}+\kappa),
\label{Eq:domain V(t) hat}
\end{equation}
with $\kappa<1/2$, then for any $n\neq\bar{n}$ we have
\begin{equation*}
2\sqrt{\hat{\eta} \hat{M}}(\bar{n}-n)\left(\sqrt{\hat{\eta} \hat{M}}(\bar{n}+n-2J)+2V(t)\right)<-\hat{C}_{\bar{n}},
\end{equation*}
where $\hat{C}_{\bar{n}}:=2\hat{\eta} \hat{M}(1-2\kappa)>0$.
Due to the continuity of the feedback controller and the fact that $u(\hat{\rho})=0$ for all $\hat{\rho}\in B_{\varepsilon}(\boldsymbol\rho_{\bar{n}})\subset B_{\xi}(\boldsymbol\rho_{\bar{n}})$ by {\textbf{H2}}, for any $V(t)$ satisfying~\eqref{Eq:domain V(t) hat}, $t\in[T_1,T_1+\delta]$ with $\delta>0$ sufficiently small and $n\neq\bar{n}$, we have $\dot{\hat{f}}_n(t)\leq -\hat{C}_{\bar{n}}\hat{f}_n(t)$. By Gr\"onwall's inequality, we have the bound  
$
\hat{f}_n(t)\leq \hat{f}_n(0)e^{-\hat{C}_{\bar{n}}t}
$
for the solution of Equation~\eqref{Eq:ODE f_hat}  starting at $\hat{f}_n(0)>0$, which implies 
\begin{equation}
\frac{1-(\hat{\rho}_v(t))_{\bar{n},\bar{n}}}{(\hat{\rho}_v(t))_{\bar{n},\bar{n}}}=\sum_{n\neq\bar{n}}\hat{f}_n(t)\leq e^{-\hat{C}_{\bar{n}}t}\sum_{n\neq\bar{n}}\hat{f}_n(0)<\frac{1-(\hat{\rho}_v(0))_{\bar{n},\bar{n}}}{(\hat{\rho}_v(0))_{\bar{n},\bar{n}}}.
\label{Eq:Estimation ODE filter}
\end{equation}
Then, $(\hat{\rho}_v(t))_{\bar{n},\bar{n}}>(\hat{\rho}_v(0))_{\bar{n},\bar{n}}$ which means that $\hat{\rho}_v(t)$ stays  in $B_{\varepsilon}(\boldsymbol\rho_{\bar{n}})$ as long as the feedback is zero. In particular a simple reasoning by contradiction shows that $u(\hat{\rho}_v(t))=0$ for all $t\geq T_1$ and, by~\eqref{Eq:Estimation ODE filter}, we have that $(\hat{\rho}_v(t))_{\bar{n},\bar{n}}$ converges to one as $t$ goes to infinity. 
Moreover, for all $t\geq T_1$ and $n\neq \bar{n}$, the dynamic of $f_n(t)$ is given by
\begin{equation*}
\dot{f}_n(t)=2\sqrt{\eta M}(\bar{n}-n)\left(\sqrt{\eta M}(\bar{n}+n-2J)+2V(t)\right)f_n(t).
\end{equation*}
By the condition~\eqref{eq:parametersrf}, there exists $V(t)$ satisfying simultaneously~\eqref{Eq:domain V(t) hat} and the following inequality
\begin{equation}
\sqrt{\eta M}(J-\bar{n}-\kappa)<V(t)<\sqrt{\eta M}(J-\bar{n}+\kappa)
\label{Eq:domain V(t)}
\end{equation}
for some $\kappa<1/2$.

Then, for any $n\neq\bar{n}$, we have 
\begin{equation*}
2\sqrt{\eta M}(\bar{n}-n)\left(\sqrt{\eta M}(\bar{n}+n-2J)+2V(t)\right)<-C_{\bar{n}}<0,
\end{equation*}
where ${C}_{\bar{n}}:=2{\eta} {M}(1-2\kappa)>0$.
Hence, we can show as above that, for all $t\geq T_1$, if $V(t)$ satisfies~\eqref{Eq:domain V(t)}--\eqref{Eq:domain V(t) hat}, then 
$
f_n(t)\leq f_n(0)e^{-C_{\bar{n}}t}
$
with $f_n(0)>0$, and $(\rho_v(t))_{\bar{n},\bar{n}}$ converges to 1 when $t$ goes to infinity. Therefore, for any $\varepsilon\in(0,\xi)$, there exists $T_2\in(0,\infty)$ and $v(t)\in\mathcal{V}$ such that, $(\rho_v(t),\hat{\rho}_v(t))$ enters $\mathbf{B}_{\varepsilon}(\boldsymbol\rho_{\bar{n}},\boldsymbol\rho_{\bar{n}})$ with $t<T_2$.
The proof is complete.\hfill$\square$
\subsubsection{Reachability of the stochastic coupled system} 
We define the stopping time $\tau^{\bar{n}}_{\varepsilon}: = \inf\{t \geq 0|\, (\rho_t,\hat{\rho}_t) \in \mathbf{B}_{\varepsilon}(\boldsymbol\rho_{\bar n},\boldsymbol\rho_{\bar n})\}$ and the compact set $\Gamma^{\bar{n}}_{\varepsilon,r}:=\left(\mathcal{S}_{N}\times\mathcal{S}_{N}\right)\setminus\left\{\mathbf{B}_{\varepsilon}(\boldsymbol\rho_{\bar{n}},\boldsymbol\rho_{\bar{n}})\cup\bigcup_{n\neq \bar{n}}\mathbf{B}_r(\boldsymbol\rho_{n},\boldsymbol\rho_{\bar{n}})\right\}$ for  $\varepsilon,r>0.$ 

Based on the results of the previous section,
we can now state the following reachability result for the stochastic coupled system~\eqref{Eq:ND SME W}--\eqref{Eq:ND SME filter W}. 

\begin{lemma}
Let $\bar{n}\in\{0,\dots,2J\}$. If $\bar{n}\in\{0,2J\}$ suppose that \emph{\textbf{H0}}, \emph{\textbf{H1}} and the condition~\eqref{eq:robustness-0,2J} hold true. Moreover assume that the hypothesis on the feedback in \hyperref[Lemma:Reachability ODE - 1]{Lemma~\ref{Lemma:Reachability ODE - 1}} is satisfied and that the initial state  $(\rho_0,\hat{\rho}_0)$ belongs to $\left(\mathcal{S}_N\times\mathcal{S}_N\right)\setminus \bigcup_{n\neq \bar{n}}(\boldsymbol\rho_{n},\boldsymbol\rho_{\bar{n}})$. Otherwise,  if $\bar{n}\in\{1,\dots,2J-1\}$, suppose that \emph{\textbf{H0}}, \emph{\textbf{H2}} and the conditions \eqref{eq:parametersr} and~\eqref{Condition u} are satisfied, and that $(\rho_0,\hat{\rho}_0)\in\mathcal{S}_N\times\mathrm{int}(\mathcal{S}_N)$.
Then, for all $\varepsilon>0$ one has 
$
\mathbb{P}_0(\tau^{\bar{n}}_{\varepsilon} < \infty)=1.
$

\label{Lemma:Reachability}
\end{lemma}
\proof
The lemma holds trivially true for $(\rho_0,\hat{\rho}_0) \in \mathbf{B}_{\varepsilon}(\boldsymbol\rho_{\bar n},\boldsymbol\rho_{\bar n})$, as in that case $\tau^{\bar{n}}_{\varepsilon} = 0$. Let us suppose $(\rho_0,\hat{\rho}_0) \in \Pi^{\bar{n}}_{\varepsilon}:=(\mathcal{S}_N\times\mathcal{S}_N)\setminus\mathbf{B}_{\varepsilon}(\boldsymbol\rho_{\bar n},\boldsymbol\rho_{\bar n})$.

Note that the condition \eqref{eq:parametersr} implies \eqref{eq:parametersrf}. Due to \hyperref[Lemma:Reachability ODE - 1]{Lemma~\ref*{Lemma:Reachability ODE - 1}}, \hyperref[Lemma:Reachability ODE - 2]{Lemma~\ref*{Lemma:Reachability ODE - 2}} and \hyperref[Thm:Support]{Theorem~\ref*{Thm:Support}}, there exist $\zeta>0$ and $T\in(0,\infty)$ such that
$
\mathbb{P}_0(\tau^{\bar{n}}_{\varepsilon}< T)\geq \zeta.
$
By the compactness of $\Gamma^{\bar{n}}_{\varepsilon,r}$ and the Feller continuity of $(\rho_t,\hat{\rho}_t)$, we have
$
\sup_{(\rho_0,\hat{\rho}_0)\in\Gamma^{\bar{n}}_{\varepsilon,r}}\mathbb{P}_0(\tau^{\bar{n}}_{\varepsilon}\geq T)\leq 1-\zeta_0<1,
$
with $\zeta\geq\zeta_0>0$. Then, by employing a similar argument as in the proof of~\cite[Theorem 4.7]{baxendale1991invariant}, we can show that there exists $K\in(0,\infty)$ such that, for all $(\rho_0,\hat{\rho}_0)\in \Gamma^{\bar{n}}_{\varepsilon,r}$ 
\begin{equation}
\mathbb{E}_0\left(\int^{\tau^{\bar{n}}_{\varepsilon}}_0 \mathds{1}_{\Gamma^{\bar{n}}_{\varepsilon,r}} (\rho_t,\hat{\rho}_t)dt\right)\leq K.
\label{Eq:Estimation Exp}
\end{equation}

Next, we construct an auxiliary function $h_{\bar{n}}(\rho,\hat{\rho})$ on the compact set $\Pi^{\bar{n}}_{\varepsilon}$ as $h_{\bar{n}}(\rho,\hat{\rho}):=-\log(1-\hat{\rho}_{\bar{n},\bar{n}})$ if $\bar{n}\in\{0,2J\}$ and \begin{equation*}
h_{\bar{n}}(\rho,\hat{\rho}):=
\begin{cases}
-\log \hat{\rho}_{n,n}, & (\rho,\hat{\rho})\in\mathbf{B}_r(\boldsymbol\rho_{n},\boldsymbol\rho_{\bar{n}})\!\cap\!(\mathcal{S}_N\times\mathrm{int}(\mathcal{S}_N)), n\neq\bar{n}\\
p_{\bar{n}}(\rho,\hat{\rho}), & (\rho,\hat{\rho})\in\Gamma^{\bar{n}}_{\varepsilon,r},
\end{cases}
\end{equation*}
if $\bar{n}\notin\{0,2J\},$ with $p_{\bar{n}}(\rho,\hat{\rho})$ any function such that $p_{\bar{n}}(\rho,\hat{\rho})\in\mathcal{C}^2(\Gamma^{\bar{n}}_{\varepsilon,r},\mathbb{R}_{\geq0}).$  Note that $h_{\bar{n}}(\rho,\hat{\rho})$ is in $\mathcal{C}^2(\Pi^{\bar{n}}_{\varepsilon}\setminus\{\bigcup_{n\neq\bar n}(\boldsymbol\rho_{n},\boldsymbol\rho_{\bar n})\},\mathbb{R}_{\geq0})$ for $\bar{n}\in\{0,2J\}$ and in $\mathcal{C}^2(\Pi^{\bar{n}}_{\varepsilon}\cap(\mathcal{S}_N\times\mathrm{int}(\mathcal{S}_N)),\mathbb{R}_{\geq0})$ for $\bar{n}\in\{1,\dots,2J-1\}.$

 Due to \hyperref[Lemma:Instability 0 and 2J]{Lemma~\ref*{Lemma:Instability 0 and 2J}} and \hyperref[Lemma:Instability General]{Lemma~\ref*{Lemma:Instability General}}, there exist two constants $r>0$ sufficiently small and $C_{\bar{n}}>0$ such that, $\Phi_{\bar{n}}(\rho,\hat{\rho}):=\mathscr{L}h_{\bar{n}}(\rho,\hat{\rho})+C_{\bar{n}}\leq 0$ for all $(\rho,\hat{\rho})\in\mathbf{B}_r(\boldsymbol\rho_{n},\boldsymbol\rho_{\bar{n}})$ and $n\neq \bar{n}$. We denote $\mathbf{C}_{\bar{n}}:=\sup_{(\rho,\hat{\rho})\in\Gamma^{\bar{n}}_{\varepsilon,r}}\Phi_{\bar{n}}(\rho,\hat{\rho})$.

For $l<r,$ let us define  $\hat\tau_{l}^{\bar n}:= \inf\{t \geq 0|\, (\rho_t,\hat{\rho}_t) \in \mathbf{B}_{l}(\boldsymbol\rho_n,\boldsymbol\rho_{\bar n})\,\textrm{for all}\, n\neq \bar{n}\}$ and denote $\sigma=t\wedge\tau^{\bar{n}}_{\varepsilon}\wedge\hat\tau_l^{\bar n}$. Under the initial conditions of the lemma, we suppose additionally that $(\rho_0,\hat{\rho}_0)\in\Gamma^{\bar{n}}_{\varepsilon,l}$.

By \hyperref[Lemma:Never Reach Lemma]{Lemma~\ref*{Lemma:Never Reach Lemma}} and \hyperref[Lemma:PosDef invariant]{Lemma~\ref*{Lemma:PosDef invariant}}, we can apply It\^o's formula on $h_{\bar{n}}(\rho_t,\hat{\rho}_t)$ for $\bar{n}\in\{0,\dots,2J\}$. By taking the expectation, we have
\begin{equation*}
\begin{split}
&\mathbb{E}_0\big(h_{\bar{n}}(\rho_{\sigma},\hat{\rho}_{\sigma})\big)-h_{\bar{n}}(\rho_{0},\hat{\rho}_{0})\\
&=\mathbb{E}_0\Big(\int^{\sigma}_0 \mathscr{L}h_{\bar{n}}(\rho_s,\hat{\rho}_s)\,ds \Big)=-C_{\bar{n}}\mathbb{E}_0(\sigma)+\mathbb{E}_0\left(\int^{\sigma}_0 \Phi_{\bar{n}}(\rho_s,\hat{\rho}_s)\,ds \right)\\
&=\!-C_{\bar{n}}\mathbb{E}_0(\sigma)\!+\!\mathbb{E}_0\Big(\int^{\sigma}_0 \mathds{1}_{\Gamma^{\bar{n}}_{\varepsilon,r}}(\rho_s,\hat{\rho}_s)\Phi_{\bar{n}}(\rho_s,\hat{\rho}_s)\,ds \Big)\!+\!\mathbb{E}_0\Big(\int^{\sigma}_0 \mathds{1}_{\Xi^{\bar{n}}_{r}}(\rho_s,\hat{\rho}_s)\Phi_{\bar{n}}(\rho_s,\hat{\rho}_s)\,ds \Big)\\
&\leq -C_{\bar{n}}\mathbb{E}_0(\sigma)+\mathbf{C}_{\bar{n}}\mathbb{E}_0\left(\int^{\sigma}_0 \mathds{1}_{\Gamma^{\bar{n}}_{\varepsilon,r}}(\rho_s,\hat{\rho}_s)\,ds \right)\leq -C_{\bar{n}}\mathbb{E}_0(\sigma)+\mathbf{C}_{\bar{n}}K,
\end{split}
\end{equation*}
where 
\begin{equation*}
\Xi^{\bar{n}}_{r}:=
\begin{cases}
\bigcup_{n\neq\bar n}\big(\mathbf{B}_r(\boldsymbol\rho_{n},\boldsymbol\rho_{\bar{n}})\setminus(\boldsymbol\rho_{n},\boldsymbol\rho_{\bar{n}})\big), &\text{if } \bar{n}\in\{0,2J\},\\
\bigcup_{n\neq \bar{n}}\mathbf{B}_r(\boldsymbol\rho_{n},\boldsymbol\rho_{\bar{n}})\cap \left(\mathcal{S}_N\times \mathrm{int}(\mathcal{S}_N)\right), &\text{if } \bar{n}\in\{1,\dots,2J-1\},
\end{cases}
\end{equation*}
and the last inequality follows from~\eqref{Eq:Estimation Exp}.

Since $h_{\bar{n}}\geq0$, $\Phi_{\bar{n}}\leq 0$ and $C_{\bar{n}}>0$, the above calculations imply 
\begin{equation}
\mathbb{E}_0(\sigma)\leq \frac{1}{C_{\bar{n}}} \big(h_{\bar{n}}(\rho_{0},\hat{\rho}_{0})+\mathbf{C}_{\bar{n}}K\big)<\infty.
\label{Eq:Estimation Exp sigma}
\end{equation}
Note that \hyperref[Lemma:Never Reach Lemma]{Lemma~\ref*{Lemma:Never Reach Lemma}} and 
\hyperref[Lemma:PosDef invariant]{Lemma~\ref*{Lemma:PosDef invariant}} imply $\mathbb{P}_0(\lim_{l\rightarrow0}\hat\tau_{l}^{\bar n}=\infty)=1$. Letting $l$ tend to zero and $t$ tend to infinity, $\sigma$ converges almost surely to $\tau^{\bar{n}}_{\varepsilon}$. By the monotone convergence theorem and the estimate~\eqref{Eq:Estimation Exp sigma}, we have
$
\mathbb{E}_0(\tau^{\bar{n}}_{\varepsilon})\leq \frac{1}{C_{\bar{n}}} \big(h_{\bar{n}}(\rho_{0},\hat{\rho}_{0})+\mathbf{C}_{\bar{n}}K\big)<\infty.
$
Then by Markov inequality, for all $(\rho_0,\hat{\rho}_0)\in(\mathcal{S}_N\times\mathcal{S}_N)\setminus \{\bigcup_{n\neq \bar{n}}(\boldsymbol\rho_{n},\boldsymbol\rho_{\bar{n}})\}$ when $\bar{n}\in\{0,2J\}$ and $(\rho_0,\hat{\rho}_0)\in\mathcal{S}_N\times\mathrm{int}(\mathcal{S}_N)$ when $\bar{n}\in\{1,\dots,2J-1\}$, we have 
\begin{equation*}
\mathbb{P}_{0}(\tau^{\bar{n}}_{\varepsilon}=\infty) = \lim_{k\rightarrow \infty} \mathbb{P}_{0}(\tau^{\bar{n}}_{\varepsilon} \geq k) \leq \lim_{k\rightarrow \infty} \mathbb{E}_{0}(\tau^{\bar{n}}_{\varepsilon})/k=0,
\end{equation*}
which implies
$
\mathbb{P}_{0}(\tau^{\bar{n}}_{\varepsilon}<\infty )=1.
$
The proof is  complete.\hfill$\square$
\subsubsection{A general result on exponential stabilization}
The following theorem provides general Lyapunov-type conditions ensuring exponential stabilization towards the target state $(\boldsymbol\rho_{\bar{n}},\boldsymbol\rho_{\bar{n}})$.
\begin{theorem}
Suppose that the assumptions of \hyperref[Lemma:Reachability]{Lemma~\ref*{Lemma:Reachability}} are satisfied. 
Additionally, assume the existence of a positive-definite function $V(\rho,\hat{\rho})$ such that $V(\rho,\hat{\rho})\!\!=\!\!0$ if and only if $(\rho,\hat{\rho})=(\boldsymbol\rho_{\bar n},\boldsymbol\rho_{\bar n})$, and $V$ is continuous on $\mathcal{S}_N\times\mathcal{S}_N$ and twice continuously differentiable on an almost surely invariant subset $\Gamma$ of $\mathcal{S}_N\times\mathcal{S}_N$ containing $\mathrm{int}(\mathcal{S}_N)\times\mathrm{int}(\mathcal{S}_N)$.
 Moreover, suppose that there exist positive constants $C$, $C_1$ and $C_2$ such that 
\begin{enumerate}
\item[(i)] $C_1 \, \mathbf{d}_B\big((\rho,\hat{\rho}),(\boldsymbol\rho_{\bar n},\boldsymbol\rho_{\bar n})\big) \leq V(\rho,\hat{\rho}) \leq C_2 \, \mathbf{d}_B\big((\rho,\hat{\rho}),(\boldsymbol\rho_{\bar n},\boldsymbol\rho_{\bar n})\big)$, for all $(\rho,\hat{\rho})\in\mathcal{S}_N\times\mathcal{S}_N$, and 
\item[(ii)] $\limsup_{(\rho,\hat{\rho})\rightarrow(\boldsymbol\rho_{\bar n},\boldsymbol\rho_{\bar n})}\frac{\mathscr{L}V(\rho,\hat{\rho})}{V(\rho,\hat{\rho})}\leq-C$.
\end{enumerate}
Then, $(\boldsymbol\rho_{\bar n},\boldsymbol\rho_{\bar n})$ is almost surely exponentially stable for the coupled system~\eqref{Eq:ND SME W}--\eqref{Eq:ND SME filter W} starting from $\Gamma$ with sample Lyapunov exponent less than or equal to $-C-\frac{K}{2}$, where $K:=\liminf_{(\rho,\hat{\rho})\rightarrow(\boldsymbol\rho_{\bar n},\boldsymbol\rho_{\bar n})}g^2(\rho,\hat{\rho})$ and $g(\rho,\hat{\rho}):=\frac{\partial V(\rho,\hat{\rho})}{\partial \rho}\frac{G_{\eta,M}(\rho)}{ V(\rho,\hat{\rho})}+\frac{\partial V(\rho,\hat{\rho})}{\partial \hat{\rho}}\frac{G_{\hat{\eta},\hat{M}}(\hat{\rho})}{ V(\rho,\hat{\rho})}$.
\label{Thm:Exp Stab General}
\end{theorem}
\textit{Sketch of the proof.} To prove \hyperref[Thm:Exp Stab General]{Theorem~\ref*{Thm:Exp Stab General}} one may follow the same steps as in~\cite[Theorem 6.2]{liang2019exponential}. In particular the presence of a function $V$ satisfying $(i)$ and such that $\mathscr{L}V\leq 0$ may be used to prove that $(\boldsymbol\rho_{\bar{n}},\boldsymbol\rho_{\bar{n}})$ is a locally stable equilibrium in probability. This, together with \hyperref[Lemma:Reachability]{Lemma~\ref*{Lemma:Reachability}} and the strong Markov property of $(\rho_t,\hat{\rho}_t)$, implies the almost sure convergence to the target equilibrium. Finally, in view of \hyperref[Lemma:PosDef invariant]{Lemma~\ref*{Lemma:PosDef invariant}}, the $\mathcal{C}^2$ regularity of the function $V$ in $\Gamma$ and the condition $(ii)$ imply 
\begin{equation*}
\limsup_{t \rightarrow \infty} \frac{1}{t} \log V(\rho_t,\hat\rho_t) \leq -C-\frac{K}2, \quad a.s.
\end{equation*}
(see~\cite[Theorem 6.2]{liang2019exponential} for more details).
The result then follows from condition~$(i)$.\qed
\subsubsection{A general result on asymptotic stabilization}
By employing similar arguments as in the first two steps of the proof in~\cite[Theorem 6.2]{liang2019exponential}, we can obtain general Lyapunov-type conditions ensuring asymptotic stabilization of the coupled system~\eqref{Eq:ND SME W}--\eqref{Eq:ND SME filter W} towards the target state. Denote $\mathcal{K}$ as the family of all continuous non-decreasing functions $\mu:\mathbb{R}_{\geq 0}\rightarrow\mathbb{R}_{\geq0}$ such that $\mu(0)=0$ and $\mu(r)>0$ for all $r>0$.
\begin{proposition}
Suppose that the assumptions of \hyperref[Lemma:Reachability]{Lemma~\ref*{Lemma:Reachability}} are satisfied.
Additionally, suppose that there exists a positive-definite function $V(\rho,\hat{\rho})$ such that $V(\rho,\hat{\rho})=0$ if and only if $(\rho,\hat{\rho})=(\boldsymbol\rho_{\bar n},\boldsymbol\rho_{\bar n})$, and $V$ is continuous on $\mathcal{S}_N\times\mathcal{S}_N$ and twice continuously differentiable on an almost surely invariant subset $\Gamma$ of $\mathcal{S}_N\times\mathcal{S}_N$ containing $\mathrm{int}(\mathcal{S}_N)\times\mathrm{int}(\mathcal{S}_N)$.
Moreover, suppose that there exists a function $\mu\in\mathcal{K}$ such that 
\begin{enumerate}
\item[(i)] $V(\rho,\hat{\rho})\geq \mu\big(\mathbf{d}_B\big((\rho,\hat{\rho}),(\boldsymbol\rho_{\bar n},\boldsymbol\rho_{\bar n})\big)\big)$, for all $(\rho,\hat{\rho})\in\mathcal{S}_N\times\mathcal{S}_N$, and 
\item[(ii)] $\mathscr{L}V(\rho,\hat{\rho})\leq0$ for all $(\rho,\hat{\rho})\in\mathbf{B}_r(\boldsymbol\rho_{\bar{n}},\boldsymbol\rho_{\bar{n}})$ with some $r>0$.
\end{enumerate}
Then, $(\boldsymbol\rho_{\bar n},\boldsymbol\rho_{\bar n})$ is almost surely asymptotically stable for the coupled system~\eqref{Eq:ND SME W}--\eqref{Eq:ND SME filter W} starting from $\Gamma.$
\label{Cor:Asym Stab General}
\end{proposition}
\begin{remark}
Following~\cite{bouten2008separation}, \hyperref[Thm:Exp Stab General]{Theorem~\ref*{Thm:Exp Stab General}} and \hyperref[Cor:Asym Stab General]{Proposition~\ref*{Cor:Asym Stab General}} can be considered as versions of the quantum separation principle in the case in which the feedback depends only on the knowledge of the estimated state. 
\end{remark}
\subsection{Explicit results on exponential stabilization and asymptotic stabilization}\label{sec:explicit}
 In this section, we establish conditions on the feedback controller $u(\hat{\rho})$ and the domain of the estimated parameters $\hat{M}$ and $\hat{\eta}$, which ensure almost sure exponential stabilization of the coupled system~\eqref{Eq:ND SME W}--\eqref{Eq:ND SME filter W} towards the target state $(\boldsymbol\rho_{\bar n},\boldsymbol\rho_{\bar n})$. We will consider separately the cases $\bar n\in\{0,2J\}$ and $\bar n\in\{1,\cdots,2J-1\}.$

\subsubsection{Stabilization results for $\bf \bar n\in\{0,2J\}$} Here, we present explicit results regarding exponential stabilization and asymptotic stabilization for the case $\bar n\in\{0,2J\}.$ 
\begin{theorem}
Consider the coupled system~\eqref{Eq:ND SME W}--\eqref{Eq:ND SME filter W} with $(\rho_{0},\hat{\rho}_0) \in (\mathcal{S}_N\times\mathcal{S}_N)\setminus\bigcup_{n\neq\bar{n}}(\boldsymbol\rho_{n},\boldsymbol\rho_{\bar{n}})$. Let $\boldsymbol \rho_{\bar{n}} \in \{\boldsymbol \rho_{0},\boldsymbol \rho_{2J}\}$ be the target state. Suppose that the assumptions on the feedback controller given in \hyperref[Lemma:Reachability ODE - 1]{Lemma~\ref*{Lemma:Reachability ODE - 1}} are satisfied, and
\begin{equation}
\frac{2N-2}{2N-1}<\sqrt{\frac{\hat{\eta}\hat{M}}{\eta M}}<\frac{1}{2}+\frac{1}{2}\sqrt{\frac{N+1}{N-1}}.
\label{Eq:CondParam0_2J}
\end{equation}
Then, $(\boldsymbol \rho_{\bar{n}},\boldsymbol \rho_{\bar{n}})$ is almost surely exponentially stable with sample Lyapunov exponent less than or equal to $-\min\{\eta M,\hat{\eta}\hat{M}\}-(N-1)\sqrt{\hat{\eta}\hat{M}}\big|\sqrt{\eta M}-\sqrt{\hat{\eta}\hat{M}}\big|$.
\label{Thm:Exp Stab 0,2J}
\end{theorem}
\proof
We consider the candidate Lyapunov function $$
V_{\bar{n}}(\rho,\hat{\rho})=\sqrt{1-\rho_{\bar{n},\bar{n}}+1-\hat{\rho}_{\bar{n},\bar{n}}}.
$$
In the following, we show that we can apply \hyperref[Thm:Exp Stab General]{Theorem~\ref*{Thm:Exp Stab General}}. We note that the set $\Gamma= (\mathcal{S}_N\times\mathcal{S}_N)\setminus\bigcup_{n\neq\bar{n}}(\boldsymbol\rho_{n},\boldsymbol\rho_{\bar{n}})$ is almost surely invariant by~\hyperref[Lemma:Never Reach Lemma]{Lemma~\ref*{Lemma:Never Reach Lemma}} and that $V_{\bar{n}}$ is twice continuously differentiable in $\Gamma$. The condition~\emph{(i)} of \hyperref[Thm:Exp Stab General]{Theorem~\ref*{Thm:Exp Stab General}} is verified because $\frac{1}{2} \, \mathbf{d}_B\big((\rho,\hat{\rho}),(\boldsymbol\rho_{\bar n},\boldsymbol\rho_{\bar n})\big) \leq V_{\bar{n}}(\rho,\hat{\rho}) \leq \, \mathbf{d}_B\big((\rho,\hat{\rho}),(\boldsymbol\rho_{\bar n},\boldsymbol\rho_{\bar n})\big)$, for all $(\rho,\hat{\rho})\in\mathcal{S}_N\times\mathcal{S}_N$.

Also, the condition~\emph{(ii)} of \hyperref[Thm:Exp Stab General]{Theorem~\ref*{Thm:Exp Stab General}} holds true, since
\begin{equation*}
\begin{split}
\mathscr{L}V_{\bar{n}}(\rho,\hat{\rho})=&\frac{1}{2V_{\bar{n}}(\rho,\hat{\rho})}\Big( u\big(\Theta_{\bar{n}}(\rho)+\Theta_{\bar{n}}(\hat{\rho})\big)-4\sqrt{\hat{\eta}\hat{M}}P_{\bar{n}}(\hat{\rho})\mathcal{T}(\rho,\hat{\rho})\hat{\rho}_{\bar{n},\bar{n}} \Big)\\
&-\frac{1}{2V^3_{\bar{n}}(\rho,\hat{\rho})}\Big(\sqrt{\hat{\eta}\hat{M}}P_{\bar{n}}(\hat{\rho})\hat{\rho}_{\bar{n},\bar{n}}+\sqrt{\eta M}P_{\bar{n}}(\rho)\rho_{\bar{n},\bar{n}}\Big)^2\\
\leq& -\mathbf{C}_{\bar{n}}(\rho,\hat{\rho})V_{\bar{n}}(\rho,\hat{\rho}),
\end{split}
\end{equation*}
with
$
\mathbf{C}_{\bar{n}}(\rho,\hat{\rho}):=-\alpha V^{m-1/2}_{\bar{n}}-2\sqrt{\hat{\eta}\hat{M}}|\mathcal{T}(\rho,\hat{\rho})|\hat{\rho}_{\bar{n},\bar{n}}+\frac{1}{2}\min\{\eta M \rho^2_{\bar{n},\bar{n}},\hat{\eta}\hat{M}\hat{\rho}^2_{\bar{n},\bar{n}}\},
$
and some constant $\alpha>0$ (see the proof of \hyperref[Lemma:Instability 0 and 2J]{Lemma~\ref*{Lemma:Instability 0 and 2J}} for the estimations of $P_{\bar{n}}$ and $\Theta_{\bar{n}}$). Thus, we have
\begin{equation}
\limsup_{(\rho,\hat{\rho})\rightarrow(\boldsymbol\rho_{\bar n},\boldsymbol\rho_{\bar n})}\frac{\mathscr{L}V_{\bar{n}}(\rho,\hat{\rho})}{V_{\bar{n}}(\rho,\hat{\rho})}\leq\limsup_{(\rho,\hat{\rho})\rightarrow(\boldsymbol\rho_{\bar n},\boldsymbol\rho_{\bar n})}-\mathbf{C}_{\bar{n}}(\rho,\hat{\rho})\leq-\bar{C}<0,
\label{Eq:Lyapunov exponent av 0_2J}
\end{equation}
where $\bar C:=\frac 12\min\{\eta M,\hat{\eta}\hat{M}\}-(N-1)\sqrt{\hat{\eta}\hat{M}}\big|\sqrt{\eta M}-\sqrt{\hat{\eta}\hat{M}}\big|.$ The positivity of $\bar C$ is guaranteed by the condition~\eqref{Eq:CondParam0_2J}. Moreover, we have the following 
\begin{equation*}
\liminf_{(\rho,\hat{\rho})\rightarrow(\boldsymbol\rho_{\bar n},\boldsymbol\rho_{\bar n})}\left(\frac{\partial V(\rho,\hat{\rho})}{\partial \rho}\frac{G_{\eta,M}(\rho)}{ V(\rho,\hat{\rho})}+\frac{\partial V(\rho,\hat{\rho})}{\partial \hat{\rho}}\frac{G_{\hat{\eta},\hat{M}}(\hat{\rho})}{ V(\rho,\hat{\rho})}\right)^2\geq \min\{\eta M,\hat{\eta}\hat{M}\}.
\end{equation*}
Then, by \hyperref[Thm:Exp Stab General]{Theorem~\ref*{Thm:Exp Stab General}}, the exponential stabilization is ensured with the exponent less than or equal to $-\bar C-\frac 12\min\{\eta M,\hat{\eta}\hat{M}\}$. The proof is then complete.\hfill$\square$

\medskip

The following result establishes  the exponential convergence towards the target state by assuming a larger  domain for the estimated parameters, compared to \hyperref[Thm:Exp Stab 0,2J]{Theorem~\ref*{Thm:Exp Stab 0,2J}}, but with  a more restrictive condition on the initial states.
\begin{theorem}
Consider the coupled system~\eqref{Eq:ND SME W}--\eqref{Eq:ND SME filter W}. Let $\boldsymbol \rho_{\bar{n}} \in \{\boldsymbol \rho_{0},\boldsymbol \rho_{2J}\}$ be the target state. Suppose that the assumptions on the feedback controller given in \hyperref[Lemma:Reachability ODE - 1]{Lemma~\ref*{Lemma:Reachability ODE - 1}} for the case $\bar{n}\in\{0,2J\}$ are satisfied, and 
\begin{equation}
\frac{2N-2}{2N-1}<\sqrt{\frac{\hat{\eta}\hat{M}}{\eta M}}<\frac{2N-2}{2N-3}.
\label{Eq:CondParamExt0_2J}
\end{equation}
Then, for all $(\rho_{0},\hat{\rho}_0) \in \mathrm{int}(\mathcal{S}_N)\times\mathrm{int}(\mathcal{S}_N)$, $(\boldsymbol \rho_{\bar{n}},\boldsymbol \rho_{\bar{n}})$ is almost surely exponentially stable with sample Lyapunov exponent less than or equal to $-\bar{C}-\bar{K}/2$ with 
\begin{equation*}
\bar{C}:=\min\big\{\frac{\eta M}{2},\frac{\hat{\eta}\hat{M}}{2}-(N-1)\sqrt{\hat{\eta}\hat{M}}\big|\sqrt{\eta M}-\sqrt{\hat{\eta}\hat{M}}\big|\big\},\quad \bar{K}:=\min\{\eta M,\hat{\eta}\hat{M}\}.
\end{equation*}
\label{Cor:Stab 0,2J}
\end{theorem}
\proof
We define
$
V_{\bar{n}}(\rho,\hat{\rho})=\sqrt{1-\rho_{\bar{n},\bar{n}}}+\sqrt{1-\hat{\rho}_{\bar{n},\bar{n}}}.
$
 Due to \hyperref[Lemma:PosDef invariant]{Lemma~\ref*{Lemma:PosDef invariant}}, $\Gamma=\mathrm{int}(\mathcal{S}_N)\times\mathrm{int}(\mathcal{S}_N)$ is almost surely invariant. Also, we note that $V_{\bar{n}}(\rho,\hat{\rho})$ is continuous on $\mathcal{S}_N\times\mathcal{S}_N$ and twice continuously differentiable on $\Gamma.$ Moreover, the condition \emph{(i)} of  \hyperref[Thm:Exp Stab General]{Theorem~\ref*{Thm:Exp Stab General}} is satisfied because we have 
$\frac{\sqrt{2}}{2}\mathbf{d}_B\big((\rho,\hat{\rho}),(\boldsymbol\rho_{\bar n},\boldsymbol\rho_{\bar n})\big) \leq V_{\bar{n}}(\rho,\hat{\rho}) \leq \mathbf{d}_B\big((\rho,\hat{\rho}),(\boldsymbol\rho_{\bar n},\boldsymbol\rho_{\bar n})\big)$ for all $(\rho,\hat{\rho})\in\mathcal{S}_N\times\mathcal{S}_N$. In addition, we show
\begin{equation*}
\begin{split}
\mathscr{L}V_{\bar{n}}(\rho,\hat{\rho})=&\frac{u\Theta_{\bar{n}}(\rho)}{2\sqrt{1-\rho_{\bar{n},\bar{n}}}}-\frac{\eta M P^2_{\bar{n}}(\rho)\rho^2_{\bar{n},\bar{n}}}{2(1-\rho_{\bar{n},\bar{n}})^{3/2}}+\frac{u\Theta_{\bar{n}}(\hat{\rho})}{2\sqrt{1-\hat{\rho}_{\bar{n},\bar{n}}}}\\
&-\frac{2\sqrt{\hat{\eta}\hat{M}}P_{\bar{n}}(\hat{\rho})\mathcal{T}(\rho,\hat{\rho})\hat{\rho}_{\bar{n},\bar{n}}}{\sqrt{1-\hat{\rho}_{\bar{n},\bar{n}}}}-\frac{\hat{\eta} \hat{M} P^2_{\bar{n}}(\hat{\rho})\hat{\rho}^2_{\bar{n},\bar{n}}}{2(1-\hat{\rho}_{\bar{n},\bar{n}})^{3/2}}\\
\leq&-\mathbf{C}_{\bar{n}}(\rho,\hat{\rho})V_{\bar{n}}(\rho,\hat{\rho})
\end{split}
\end{equation*}
where 
\begin{equation*}
\mathbf{C}_{\bar{n}}(\rho,\hat{\rho}):=\min\big\{\frac{\eta M \rho^2_{\bar{n},\bar{n}}}{2}, \frac{\hat{\eta} \hat{M} \hat{\rho}^2_{\bar{n},\bar{n}}}{2}-2|\mathcal{T}(\rho,\hat{\rho})| \sqrt{\hat{\eta}\hat{M}}\hat{\rho}_{\bar{n},\bar{n}}+\alpha(1-\hat{\rho}_{\bar{n},\bar{n}})^{m-\frac12}\big\},
\end{equation*}
with some constant $\alpha>0,$ which can be determined by the estimations provided in the proof of \hyperref[Lemma:Instability 0 and 2J]{Lemma~\ref*{Lemma:Instability 0 and 2J}}. Thus, we have
\begin{equation*}
\limsup_{(\rho,\hat{\rho})\rightarrow(\boldsymbol\rho_{\bar n},\boldsymbol\rho_{\bar n})}\frac{\mathscr{L}V_{\bar{n}}(\rho,\hat{\rho})}{V_{\bar{n}}(\rho,\hat{\rho})}\leq\limsup_{(\rho,\hat{\rho})\rightarrow(\boldsymbol\rho_{\bar n},\boldsymbol\rho_{\bar n})}-\mathbf{C}_{\bar{n}}(\rho,\hat{\rho})\leq-\bar{C}<0,
\end{equation*}
where $\bar{C}:=\min\big\{\frac{\eta M}{2}, \frac{\hat{\eta} \hat{M}}{2}-2J\sqrt{\hat{\eta}\hat{M}}\big|\sqrt{\eta M}-\sqrt{\hat{\eta}\hat{M}}\big|\big\}>0$. The positivity of $\bar C$ is guaranteed by the condition~\eqref{Eq:CondParamExt0_2J}. Hence, the condition \emph{(ii)} of \hyperref[Thm:Exp Stab General]{Theorem~\ref*{Thm:Exp Stab General}} is satisfied. As a consequence \hyperref[Thm:Exp Stab General]{Theorem~\ref*{Thm:Exp Stab General}} can be applied and, in order to find the sample Lyapunov exponent, we notice that
\begin{equation*}
\liminf_{(\rho,\hat{\rho})\rightarrow(\boldsymbol\rho_{\bar n},\boldsymbol\rho_{\bar n})}\left(\frac{\partial V(\rho,\hat{\rho})}{\partial \rho}\frac{G_{\eta,M}(\rho)}{ V(\rho,\hat{\rho})}+\frac{\partial V(\rho,\hat{\rho})}{\partial \hat{\rho}}\frac{G_{\hat{\eta},\hat{M}}(\hat{\rho})}{ V(\rho,\hat{\rho})}\right)^2\geq\bar{K}.
\end{equation*}
The proof is complete.\hfill$\square$

\medskip
In the following, with a less restrictive assumption on the initial condition, we show the asymptotic stabilization of the target state. 
\begin{proposition}
Consider the coupled system~\eqref{Eq:ND SME W}--\eqref{Eq:ND SME filter W} with  $(\rho_{0},\hat{\rho}_0) \in \mathcal{S}_N\times(\mathcal{S}_N\setminus\{\boldsymbol\rho_{\bar{n}}\})$. Suppose that the assumptions on the feedback controller  in \hyperref[Lemma:Reachability ODE - 1]{Lemma~\ref*{Lemma:Reachability ODE - 1}} and the condition~\eqref{Eq:CondParam0_2J} are satisfied. Then, $(\boldsymbol \rho_{\bar{n}},\boldsymbol \rho_{\bar{n}})$ is almost surely asymptotically stable.
\end{proposition}
\proof
It is sufficient to consider $\mathsf{V}_{\bar{n}}(\rho,\hat{\rho})=1-\rho_{\bar{n},\bar{n}}+\sqrt{1-\hat{\rho}_{\bar{n},\bar{n}}}$. Due to \hyperref[Lemma:Never Reach Lemma]{Lemma~\ref*{Lemma:Never Reach Lemma}}, $\Gamma=\mathcal{S}_N\times(\mathcal{S}_N\setminus\{\boldsymbol\rho_{\bar{n}}\})$ is almost surely invariant. Moreover, the function $\mathsf{V}_{\bar{n}}(\rho,\hat{\rho})$ is continuous on $\mathcal{S}_N\times\mathcal{S}_N$ and twice continuously differentiable on $\Gamma.$ 
Also, the function V satisfies 
$
\mathsf{V}_{\bar{n}}(\rho,\hat{\rho}) \geq \mathbf{d}^2_B\big((\rho,\hat{\rho}),(\boldsymbol\rho_{\bar n},\boldsymbol\rho_{\bar n})\big)
$
for all $(\rho,\hat{\rho})\in\mathcal{S}_N\times\mathcal{S}_N$ and 
\begin{equation*}
\begin{split}
\mathscr{L}\mathsf{V}_{\bar{n}}(\rho,\hat{\rho})=&u\Theta_{\bar{n}}(\rho)+\frac{u\Theta_{\bar{n}}(\hat{\rho})}{2\sqrt{1-\hat{\rho}_{\bar{n},\bar{n}}}}-\frac{2\sqrt{\hat{\eta}\hat{M}}P_{\bar{n}}(\hat{\rho})\mathcal{T}(\rho,\hat\rho)\hat{\rho}_{\bar{n},\bar{n}}}{\sqrt{1-\hat{\rho}_{\bar{n},\bar{n}}}}-\frac{\hat{\eta} \hat{M} P^2_{\bar{n}}(\hat{\rho})\hat{\rho}^2_{\bar{n},\bar{n}}}{2(1-\hat{\rho}_{\bar{n},\bar{n}})^{3/2}}\\
\leq&-\mathsf{C}_{\bar{n}}(\rho,\hat{\rho})\sqrt{1-\hat{\rho}_{\bar{n},\bar{n}}}
\end{split}
\end{equation*}
where 
\begin{equation*}
\mathsf{C}_{\bar{n}}(\rho,\hat{\rho}):=\frac{\hat{\eta} \hat{M} \hat{\rho}^2_{\bar{n},\bar{n}}}{2}-2\sqrt{\hat{\eta}\hat{M}}\hat{\rho}_{\bar{n},\bar{n}}|\mathcal{T}(\rho,\hat{\rho})| +\alpha(1-\hat{\rho}_{\bar{n},\bar{n}})^{m-\frac12}\big(1+\sqrt{1-\rho_{\bar{n},\bar{n}}}\big).
\end{equation*}
Moreover, we have
\begin{equation*}
\begin{split}
\liminf_{(\rho,\hat{\rho})\rightarrow(\boldsymbol\rho_{\bar n},\boldsymbol\rho_{\bar n})}\mathsf{C}_{\bar{n}}(\rho,\hat{\rho})\geq \frac{\hat{\eta} \hat{M}}{2}-2J\sqrt{\hat{\eta}\hat{M}}\left|\sqrt{\eta M}-\sqrt{\hat{\eta}\hat{M}}\right|>0,
\end{split}
\end{equation*}
the positivity of the last term is guaranteed by the condition~\eqref{Eq:CondParam0_2J}. Then, there always exists $r>0$ such that $\mathscr{L}\mathsf{V}_{\bar{n}}(\rho,\hat{\rho})\leq0$ for all $(\rho,\hat{\rho})\in\mathbf{B}_r(\boldsymbol\rho_{\bar{n}},\boldsymbol\rho_{\bar{n}})$. Thus we can apply
\hyperref[Cor:Asym Stab General]{Proposition~\ref*{Cor:Asym Stab General}} to conclude the proof.\hfill$\square$
\subsubsection{Stabilization results for $\bf\bar{n}\in\{1,\dots,2J-1\}$} Now, we present explicit results regarding exponential sabilization and asymptotic stabilization for the case $\bar{n}\in\{1,\dots,2J-1\}$. We define 
$C_{\bar{n}}:=\min\big\{\frac{\eta M}{2}, \frac{\hat{\eta} \hat{M}}{2}-\frac{\sqrt{\hat{\eta}\hat{M}}L_{\bar{n}}}{2}\big|\sqrt{\eta M}-\sqrt{\hat{\eta}\hat{M}}\big|\big\}$, with $L_{\bar{n}}:=4|J-\bar{n}|\max\{\bar{n},2J-\bar{n}\}.$
\begin{theorem}
Consider the coupled system~\eqref{Eq:ND SME W}--\eqref{Eq:ND SME filter W} with $(\rho_{0},\hat{\rho}_0) \in \mathrm{int}(\mathcal{S}_N)\times\mathrm{int}(\mathcal{S}_N)$. Let $\boldsymbol \rho_{\bar{n}} \in \{\boldsymbol \rho_{1},\dots,\boldsymbol \rho_{2J-1}\}$ be the target state. Suppose that the assumptions on the feedback controller given in \hyperref[Lemma:Reachability ODE - 2]{Lemma~\ref*{Lemma:Reachability ODE - 2}} are satisfied, and 
\begin{equation}
\begin{cases}
(N-1)\sqrt{\eta M}>(N-2)\sqrt{\hat{\eta}\hat{M}}>(N-3)\sqrt{\eta M}&\text{ if }\bar{n}=J;\\
\frac{L_{\bar{n}}}{L_{\bar{n}}-1}\sqrt{\eta M}>\sqrt{\hat{\eta}\hat{M}}>\frac{L_{\bar{n}}}{L_{\bar{n}}+1}\sqrt{\eta M}&\text{ if }\bar{n}\neq J.
\end{cases}
\label{Eq:Cond Param 1...2J-1}
\end{equation}
Then, $(\boldsymbol \rho_{\bar{n}},\boldsymbol \rho_{\bar{n}})$ is almost surely exponentially stable with sample Lyapunov exponent less than or equal to $-C_{\bar{n}}.$
\label{Thm:Exp Stab 1...2J-1}
\end{theorem}
\proof
Consider the following candidate Lyapunov function
\begin{equation*}
\mathbf{V}_{\bar{n}}(\rho,\hat{\rho})=V_{\bar{n}}(\rho)+V_{\bar{n}}(\hat{\rho})=\sum_{n\neq\bar{n}}\sqrt{\rho_{n,n}}+\sum_{n\neq\bar{n}}\sqrt{\hat{\rho}_{n,n}}.
\end{equation*}
Due to \hyperref[Lemma:PosDef invariant]{Lemma~\ref*{Lemma:PosDef invariant}}, $\Gamma=\mathrm{int}(\mathcal{S}_N)\times\mathrm{int}(\mathcal{S}_N)$ is almost surely invariant. The function $\mathbf{V}_{\bar{n}}(\rho,\hat{\rho})$ is continuous on $\mathcal{S}_N\times\mathcal{S}_N$ and twice continuously differentiable on $\Gamma.$

By applying Jensen inequality, for all $(\rho,\hat{\rho})\in\mathcal{S}_N\times\mathcal{S}_N,$ we can show

$
\frac{\sqrt{2}}{2} \, \mathbf{d}_B\big((\rho,\hat{\rho}),(\boldsymbol\rho_{\bar n},\boldsymbol\rho_{\bar n})\big) \leq \mathbf{V}_{\bar{n}}(\rho,\hat{\rho}) \leq \sqrt{2J} \, \mathbf{d}_B\big((\rho,\hat{\rho}),(\boldsymbol\rho_{\bar n},\boldsymbol\rho_{\bar n})\big).
$

Based on the estimates on $u$, $\Theta_{\bar{n}}$ and $P_{\bar{n}}$ in the proof of \hyperref[Lemma:Instability General]{Lemma~\ref*{Lemma:Instability General}}, and the fact that $u=0$ for all $(\rho,\hat{\rho})$ in a sufficiently small neighbourhood of the target state (hypothesis \textbf{H2}), we have the following estimate on the infinitesimal generator of $\mathbf{V}_{\bar{n}}(\rho,\hat{\rho})$ for all $(\rho,\hat{\rho})\in\mathbf{B}_r(\boldsymbol\rho_{\bar{n}},\boldsymbol\rho_{\bar{n}})$ with $r>0$ sufficiently small,
\begin{equation*}
\begin{split}
\mathscr{L}\mathbf{V}_{\bar{n}}(\rho,\hat{\rho})\leq&-\frac{\eta M}{2}\big(1-|P_{\bar{n}}(\rho)|)^2V_{\bar{n}}(\rho)-\frac{\hat{\eta} \hat{M}}{2}\big(1-|P_{\bar{n}}(\hat{\rho})|)^2V_{\bar{n}}(\hat{\rho})\\
&+2\sqrt{\hat{\eta} \hat{M}}\big(l_{\bar{n}}+|P_{\bar{n}}(\hat{\rho})|\big)|\mathcal{T}(\rho,\hat{\rho})|V_{\bar{n}}(\hat{\rho})\\
\leq&-\mathbf{C}_{\bar{n}}(\rho,\hat{\rho})\mathbf{V}_{\bar{n}}(\rho,\hat{\rho})
\end{split}
\end{equation*}
where $l_{\bar{n}}:=\max\{\bar{n},2J-\bar{n}\}$ and 
\begin{align}
\mathbf{C}_{\bar{n}}(\rho,\hat{\rho}):=\min&\left\lbrace\frac{\eta M}{2}\big(1-|P_{\bar{n}}(\rho)|)^2, \right.\label{Eq:C 1...2J-1}\\
&\quad\left.\frac{\hat{\eta} \hat{M}}{2}\big(1-|P_{\bar{n}}(\hat{\rho})|)^2-2\sqrt{\hat{\eta} \hat{M}}\big(l_{\bar{n}}+|P_{\bar{n}}(\hat{\rho})|\big)|\mathcal{T}(\rho,\hat{\rho})|\right\rbrace.\nonumber
\end{align}
Thus, we have
\begin{equation}
\limsup_{(\rho,\hat{\rho})\rightarrow(\boldsymbol\rho_{\bar n},\boldsymbol\rho_{\bar n})}\frac{\mathscr{L}\mathbf{V}_{\bar{n}}(\rho,\hat{\rho})}{\mathbf{V}_{\bar{n}}(\rho,\hat{\rho})}\leq\limsup_{(\rho,\hat{\rho})\rightarrow(\boldsymbol\rho_{\bar n},\boldsymbol\rho_{\bar n})}-\mathbf{C}_{\bar{n}}(\rho,\hat{\rho})\leq-C_{\bar{n}}<0,
\label{Eq:Lyapunov exponent av 1...2J-1}
\end{equation}
where the positivity of $C_{\bar{n}}$ is guaranteed by the condition~\eqref{Eq:Cond Param 1...2J-1}. Thus we can apply \hyperref[Thm:Exp Stab General]{Theorem~\ref*{Thm:Exp Stab General}} and the proof is complete.\hfill$\square$

\medskip

In the following, we show the asymptotic stabilization of the target state with a weaker condition on the initial states. 
\begin{proposition}
Consider the coupled system~\eqref{Eq:ND SME W}--\eqref{Eq:ND SME filter W} with $(\rho_{0},\hat{\rho}_0) \in \mathcal{S}_N\times\mathrm{int}(\mathcal{S}_N)$. Let $\boldsymbol \rho_{\bar{n}} \in \{\boldsymbol \rho_{1},\dots,\boldsymbol \rho_{2J-1}\}$ be the target state. Suppose that the assumptions on the feedback controller in \hyperref[Lemma:Reachability ODE - 2]{Lemma~\ref*{Lemma:Reachability ODE - 2}} for $\bar{n}\in\{1,\dots,2J-1\}$ and the condition~\eqref{Eq:Cond Param 1...2J-1} are satisfied. Then, $(\boldsymbol \rho_{\bar{n}},\boldsymbol \rho_{\bar{n}})$ is almost surely asymptotically stable.
\label{Cor:Asy Stab 1...2J-1}
\end{proposition}
\proof
Consider the following candidate Lyapunov function
\begin{equation*}
\mathbf{V}_{\bar{n}}(\rho,\hat{\rho})=1-\rho_{\bar{n},\bar{n}}+\sum_{n\neq\bar{n}}\sqrt{\hat{\rho}_{n,n}}.
\end{equation*}
Due to \hyperref[Lemma:PosDef invariant]{Lemma~\ref*{Lemma:PosDef invariant}}, $\mathcal{S}_N\times\mathrm{int}(\mathcal{S}_N)$ is almost surely invariant. The function $\mathbf{V}_{\bar{n}}(\rho,\hat{\rho})$ is continuous on $\mathcal{S}_N\times\mathcal{S}_N$ and twice continuously differentiable on $\mathcal{S}_N\times\mathrm{int}(\mathcal{S}_N)$. The result can be shown by applying \hyperref[Cor:Asym Stab General]{Proposition~\ref*{Cor:Asym Stab General}}, since by applying Jensen inequality, we can easily show  
$
\mathbf{V}_{\bar{n}}(\rho,\hat{\rho}) \geq \mathbf{d}^2_B\big((\rho,\hat{\rho}),(\boldsymbol\rho_{\bar n},\boldsymbol\rho_{\bar n})\big),
$
for all $(\rho,\hat{\rho})\in\mathcal{S}_N\times\mathcal{S}_N.$
Also, based on the same calculations provided in the proof of \hyperref[Thm:Exp Stab 1...2J-1]{Theorem~\ref*{Thm:Exp Stab 1...2J-1}}, for all $(\rho,\hat{\rho})\in\mathbf{B}_r(\boldsymbol\rho_{\bar{n}},\boldsymbol\rho_{\bar{n}})$ with $r>0$ sufficiently small, we have 
$$
\mathscr{L}\mathbf{V}_{\bar{n}}(\rho,\hat{\rho})\leq-\mathbf{C}_{\bar{n}}(\rho,\hat{\rho})\sum_{n\neq\bar{n}}\sqrt{\hat{\rho}_{n,n}},
$$
where $\mathbf{C}_{\bar{n}}(\rho,\hat{\rho})$ is defined in~\eqref{Eq:C 1...2J-1}. Moreover
$$
\liminf_{(\rho,\hat{\rho})\rightarrow(\boldsymbol\rho_{\bar n},\boldsymbol\rho_{\bar n})}\mathbf{C}_{\bar{n}}(\rho,\hat{\rho})\geq C_{\bar{n}}>0.
$$
Thus, there always exists $r>0$ such that, $\mathscr{L}\mathbf{V}_{\bar{n}}(\rho,\hat{\rho})\leq0$ for all $(\rho,\hat{\rho})\in\mathbf{B}_r(\boldsymbol\rho_{\bar{n}},\boldsymbol\rho_{\bar{n}})$. Then the proof is complete.\hfill$\square$
\begin{remark}
Under the assumptions of the above proposition, the exponential convergence can be ensured by adapting the construction of the Lyapunov function proposed in~\cite{cardona2020exponential} to our case. However, obtaining an estimate of the convergence rate with this method appears to be difficult. 
\end{remark}
\subsection{Parametrized feedback laws} \label{sec:par} As an example of application of the previous results, we design parametrized feedback laws which stabilize exponentially $(\rho_t,\hat{\rho}_t)$ almost surely towards some predetermined target eigenstate $(\boldsymbol\rho_{\bar n},\boldsymbol\rho_{\bar n})$.

We start by considering the special case $\bar n\in\{0,2J\}$.
\begin{theorem}
Consider the coupled system~\eqref{Eq:ND SME W}--\eqref{Eq:ND SME filter W} with $(\rho_{0},\hat{\rho}_0) \in (\mathcal{S}_N\times\mathcal{S}_N)\setminus\bigcup_{n\neq\bar{n}}(\boldsymbol\rho_{n},\boldsymbol\rho_{\bar{n}})$. Let $\boldsymbol \rho_{\bar{n}} \in \{\boldsymbol \rho_{0},\boldsymbol \rho_{2J}\}$ be the target state. Suppose that the condition~\eqref{Eq:CondParam0_2J} is satisfied, and define
the feedback controller
\begin{equation}
u_{\bar{n}}(\hat{\rho})  = \alpha \big(1-\mathrm{Tr}(\hat{\rho} \boldsymbol \rho_{\bar{n}})\big)^{\beta},
\label{Eq:u 0_2J}
\end{equation}
where $\alpha>0$ and $\beta \geq1$. Then, $(\boldsymbol \rho_{\bar{n}},\boldsymbol \rho_{\bar{n}})$ is almost surely exponentially stable with sample Lyapunov exponent less than or equal to the value defined in \hyperref[Thm:Exp Stab 0,2J]{Theorem~\ref*{Thm:Exp Stab 0,2J}}. 
\label{Thm:Feedback Exp Stab 0,2J}
\end{theorem}
We remark that \hyperref[Thm:Feedback Exp Stab 0,2J]{Theorem~\ref*{Thm:Feedback Exp Stab 0,2J}} provides a proof of~\cite[Conjecture~4.4]{liang2019exponential_TwoUnknown} in the case $\boldsymbol \rho_{\bar{n}} \in \{\boldsymbol \rho_{0},\boldsymbol \rho_{2J}\}$. Note that, unlike the present paper, in~\cite{liang2019exponential_TwoUnknown} the physical parameters were supposed to be known.

To tackle the case in which $\bar n$ is not necessarily equal to $0$ or $2J$, and in order to construct a feedback controller satisfying \textbf{H2}, we define a continuously differentiable function $f:[0,1]\to[0,1]$ as follows
\begin{equation*}
f(x) = 
\begin{cases}
0,&\text{if }x\in[0,\epsilon_1);\\
\frac12\sin\left(\frac{\pi(2x-\epsilon_1-\epsilon_2)}{2(\epsilon_2-\epsilon_1)}\right)+\frac12,&\text{if }x\in[\epsilon_1,\epsilon_2);\\
1,&\text{if }x\in(\epsilon_2,1],
\end{cases}
\end{equation*}
where $0<\epsilon_1<\epsilon_2<1$. We then have the following result.
\begin{theorem}
Consider the coupled system~\eqref{Eq:ND SME W}--\eqref{Eq:ND SME filter W} with $(\rho_{0},\hat{\rho}_0) \in \mathrm{int}(\mathcal{S}_N)\times\mathrm{int}(\mathcal{S}_N)$. Let $\boldsymbol \rho_{\bar{n}} \in \{\boldsymbol \rho_{0},\dots,\boldsymbol \rho_{2J}\}$ be the target state. Suppose that the condition~\eqref{Eq:CondParam0_2J} is satisfied for $\bar{n}\in\{0,2J\}$ and the condition~\eqref{Eq:Cond Param 1...2J-1} is satisfied for $\bar{n}\in\{1,\dots,2J-1\}$. Define the feedback controller
\begin{equation}
u_{\bar{n}}=\alpha \big(J-\bar{n}-\mathrm{Tr}(J_z\hat{\rho})\big)^{\beta}f(1-\hat{\rho}_{\bar{n},\bar{n}}),
\label{Eq:u 1...2J-1}
\end{equation}
where $\alpha>0$ and $\beta\geq1$. Then, $(\boldsymbol \rho_{\bar{n}},\boldsymbol \rho_{\bar{n}})$ is almost surely exponentially stable with sample Lyapunov exponent less than or equal to the value defined in \hyperref[Cor:Stab 0,2J]{Theorem~\ref*{Cor:Stab 0,2J}} for $\bar{n}\in\{0,2J\}$ and in \hyperref[Thm:Exp Stab 1...2J-1]{Theorem~\ref*{Thm:Exp Stab 1...2J-1}} for $\bar{n}\in\{1,\dots,2J-1\}$. 
\label{Thm: Feedback Exp Stab 0...2J}
\end{theorem}
\begin{remark}
If the conjecture proposed in~\cite[Remark 6.6]{liang2019exponential} holds true, then the initial condition  $(\rho_{0},\hat{\rho}_0)$ in \hyperref[Thm: Feedback Exp Stab 0...2J]{Theorem~\ref*{Thm: Feedback Exp Stab 0...2J}} can be taken in $(\mathcal{S}_N\times\mathcal{S}_N)\setminus\bigcup_{n\neq\bar{n}}(\boldsymbol\rho_{n},\boldsymbol\rho_{\bar{n}})$. This assertion has been verified for the two-level case in~\cite{liang2020robustness_two}.
\end{remark}
\section{Simulations}\label{sec:sim}
\begin{figure}[thpb]
\centering
\includegraphics[width=12.5cm]{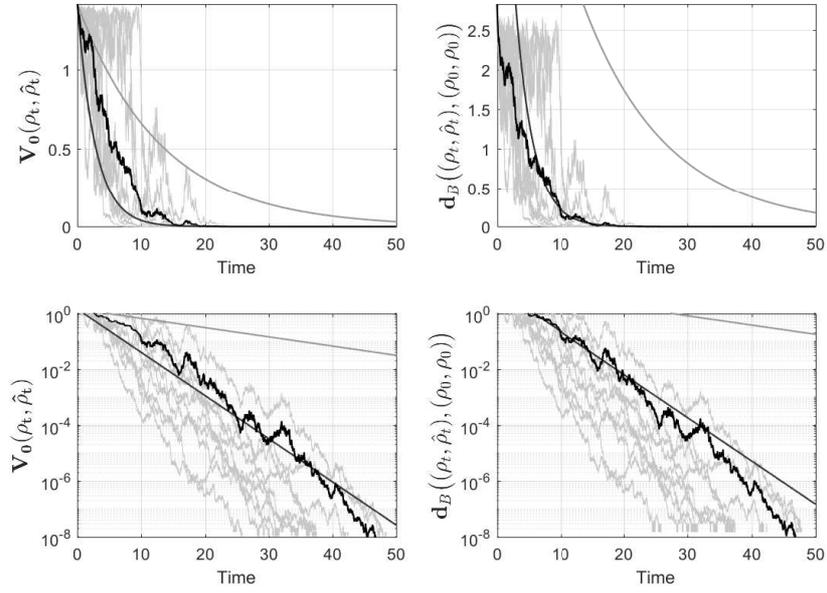}
\caption{Exponential stabilization of the coupled 3-level quantum angular momentum system towards $(\boldsymbol \rho_0,\boldsymbol \rho_0)$ with the feedback controller~\eqref{Eq:u 0_2J} starting at $(\rho_0,\hat{\rho}_0)=(\boldsymbol \rho_2,\boldsymbol \rho_1)$ with $\omega=0.4$, $\eta=0.4$, $M=1.4$, $\hat{\omega}=0.5$, $\hat{\eta}=0.5$, $\hat{M}=1.5$, $\alpha = 5$ and $\beta = 2$: the black curves represent the mean value of 10 arbitrary sample trajectories, the light grey curves represent the exponential reference with exponent $\nu_{\textrm{av}}=-0.0761$, the dark grey curve represents the exponential reference with exponent $\nu_{\textrm{s}}=-0.3561$. The figures at the bottom are the semi-log versions of the ones at the top.}
\label{Fig:3DRobust_0}
\end{figure}

\begin{figure}[thpb]
\centering
\includegraphics[width=12.5cm]{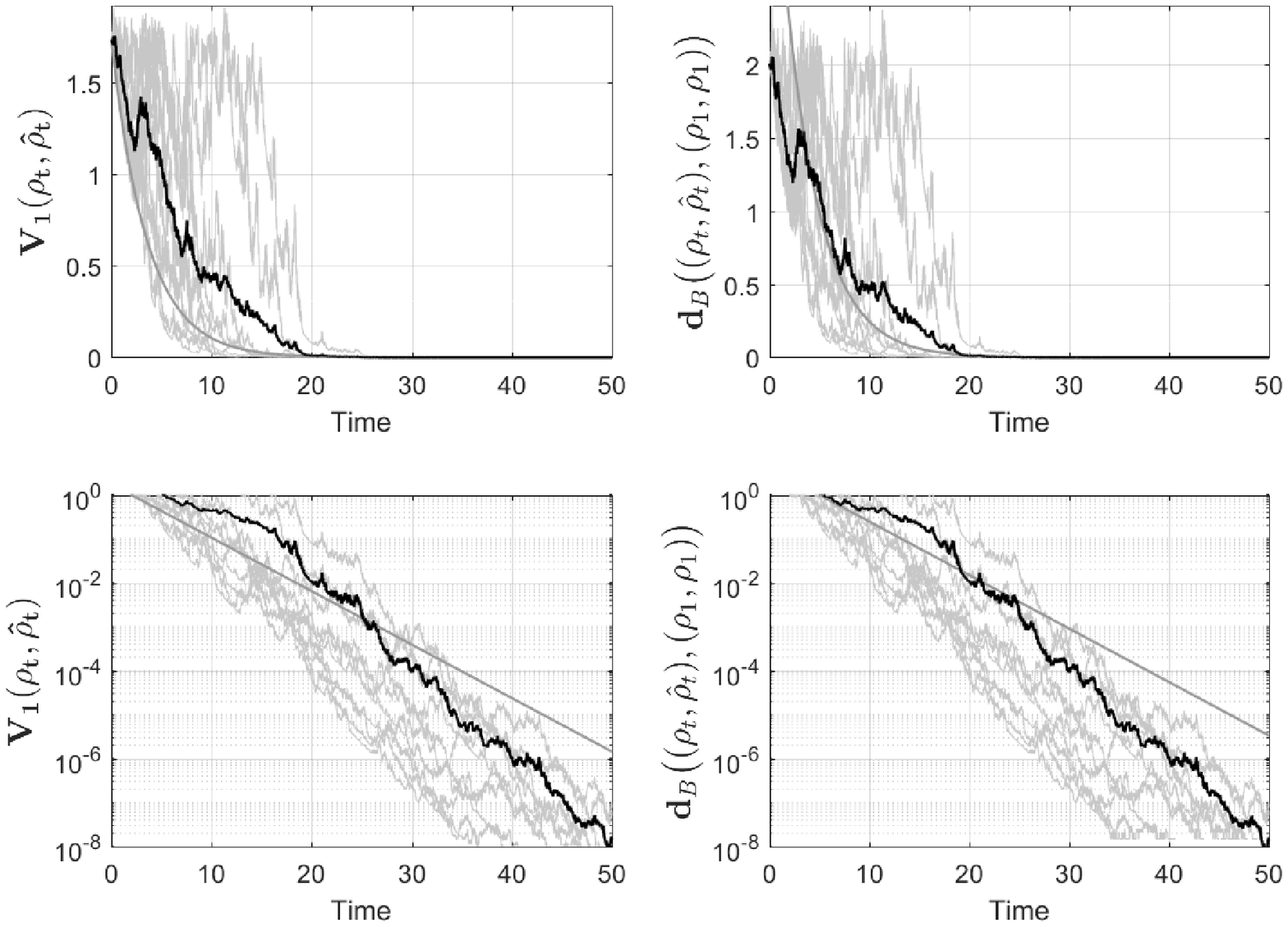}
\caption{Exponential stabilization of the coupled 3-level quantum angular momentum system towards $(\boldsymbol \rho_1,\boldsymbol \rho_1)$ with the feedback controller~\eqref{Eq:u 1...2J-1} starting at $(\rho_0,\hat{\rho}_0)=(\mathrm{diag}(0.2,0.2,0.6),\mathrm{diag}(0.3,0.3,0.4))\in\mathrm{int}(\mathcal{S}_N)\times\mathrm{int}(\mathcal{S}_N)$ with $\omega=0.4$, $\eta=0.4$, $M=1.4$, $\hat{\omega}=0.5$, $\hat{\eta}=0.5$, $\hat{M}=1.5$, $\alpha = 5$ and $\beta = 2$: the black curves represent the mean value of 10 arbitrary sample trajectories, the light grey curves represent the exponential reference with exponent $\nu_{\textrm{s}}=-0.28$. The figures at the bottom are the semi-log versions of the ones at the top.}
\label{Fig:3DRobust_1}
\end{figure}

In this section, we illustrate our results by numerical simulations in the case of a coupled three-level quantum angular momentum system. In this case $N=3$ and $J=1$. The values of the physical and experimental parameters are chosen as $\omega=0.4$, $\eta=0.4$, $M=1.4$, $\hat{\omega}=0.5$, $\hat{\eta}=0.5$, $\hat{M}=1.5$. We first illustrate the convergence of the coupled system~\eqref{Eq:ND SME W}--\eqref{Eq:ND SME filter W}  starting at $(\rho_0,\hat{\rho}_0)=(\boldsymbol \rho_2,\boldsymbol \rho_1)$ towards the target state $(\boldsymbol\rho_0,\boldsymbol\rho_0)$ by applying a feedback controller of the form~\eqref{Eq:u 0_2J}, with $\alpha = 5$ and $\beta = 2$. This is shown in \hyperref[Fig:3DRobust_0]{Figure~\ref*{Fig:3DRobust_0}}. Then, in \hyperref[Fig:3DRobust_1]{Figure~\ref*{Fig:3DRobust_1}}, we show the convergence of the coupled system starting at $(\rho_0,\hat{\rho}_0)=(\mathrm{diag}(0.2,0.2,0.6),\mathrm{diag}(0.3,0.3,0.4))\in\mathrm{int}(\mathcal{S}_N)\times\mathrm{int}(\mathcal{S}_N)$, towards the target state $(\boldsymbol\rho_1,\boldsymbol\rho_1)$ by a feedback controller of the form~\eqref{Eq:u 1...2J-1}, with $\alpha = 5$ and $\beta = 2$.

By Equation~\eqref{Eq:Lyapunov exponent av 0_2J} and Equation~\eqref{Eq:Lyapunov exponent av 1...2J-1}, heuristically we have that the rate of convergence of the expectation of the Lyapunov function is less than or equal to $\nu_{\textrm{av}}=-\frac{1}{2}\min\{\eta M ,\hat{\eta}\hat{M}\}+2\sqrt{\hat{\eta}\hat{M}}\big|\sqrt{\eta M}-\sqrt{\hat{\eta}\hat{M}}\big|$ for $\bar{n}\in\{0,2\}$, and $\nu_{\textrm{av}}=-\frac{1}{2}\min\{\eta M, \hat{\eta} \hat{M}\}$ for $\bar{n}=1$. This property is confirmed through simulations, see \hyperref[Fig:3DRobust_0]{Figure~\ref*{Fig:3DRobust_0}} and \hyperref[Fig:3DRobust_1]{Figure~\ref*{Fig:3DRobust_1}}. In the figures, the light grey curves represent the exponential reference with the exponent $\nu_{\textrm{av}}$ and the black curves describe the mean values of the Lyapunov functions (Bures distances) of ten samples. In the figures, in particular in the semi-log versions, we can see that the black and the light grey curves have similar asymptotic behaviors. In \hyperref[Fig:3DRobust_0]{Figure~\ref*{Fig:3DRobust_0}}, we observe that the dark grey curves describing the exponential reference with exponent $\nu_{\textrm{s}}:=\nu_{\textrm{av}}-\frac{1}{2}\min\{\eta M, \hat{\eta} \hat{M}\}$ have similar asymptotic behaviors compared to ten sample trajectories. Note that for $\bar{n}=1,$ the same Lyapunov exponent $\nu_{\textrm{s}}$  obtained from \hyperref[Thm:Exp Stab 1...2J-1]{Theorem~\ref*{Thm:Exp Stab 1...2J-1}} coincides with $\nu_{\textrm{av}}.$
\section{Conclusion and perspectives}
In this paper, we proved a general exponential stabilization result for $N$-level quantum angular momentum systems, robust with respect to imprecise choices of the estimated physical parameters and wrong initialization of the estimated state. Such a robustness property was obtained by analyzing the asymptotic behavior of the coupled system describing the evolution of the quantum state and the associated estimated state. More precisely, we showed exponential stabilization of the coupled system towards a pair $(\boldsymbol \rho_{\bar{n}},\boldsymbol \rho_{\bar{n}})$, with $\boldsymbol \rho_{\bar{n}}$ being a chosen eigenstate of the measurement operator $J_z$. 
Furthermore, we gave explicit examples of stabilizing feedback control laws. 
Future research lines will concern the robustness properties of the stabilizing feedback controller in presence of delays for $N$-level quantum angular momentum systems and the adaptation of the robust exponential stabilization results to general open quantum systems.


\bibliographystyle{siamplain}
\bibliography{Ref_Thesis_LIANG}
\end{document}